\documentclass[11pt]{article}
\usepackage{geometry}
\usepackage{graphicx}
\usepackage{amsmath}
\usepackage{amssymb}
\usepackage{amsthm}
\usepackage{float}
\usepackage{epstopdf}
\usepackage{color}
\usepackage{colordvi}
\usepackage{algorithm}
\usepackage{algpseudocode}
\usepackage{subcaption}
\usepackage{siunitx}
\usepackage{verbatim}

\usepackage{tikz}
\usetikzlibrary{angles,quotes}
\usetikzlibrary{calc}
\usetikzlibrary{babel}
\graphicspath{{Plots/SimpleCrack/}{Plots/WheelRail/}{Plots/}}

\newcommand{\R}{\mathbb{R}}

\renewcommand{\Vec}[1]{{\boldsymbol{#1}}}
\newcommand{\Mat}[1]{{\boldsymbol{#1}}}

\DeclareMathOperator*{\argmin}{arg\,min} 

\definecolor{magenta}{rgb}{.5,0,.5}
\definecolor{black}{rgb}{1.0,1.0,1.0}
\definecolor{magenta}{rgb}{.1,0,.3}
\definecolor{gruen}{rgb}{0.2,0.5,.5}
\definecolor{light}{rgb}{ 0.992, 0.961,  0.902}
\definecolor{Tan}{rgb}{ 0.992, 0.9,  0.902}

\title{A Physics-Based Model Reduction Approach for Node-to-Segment Contact Problems in Linear Elasticity}

\author{Diana Manvelyan\footnote{Siemens AG, Technology, Munich , Germany, diana.manvelyan@siemens.com},    
Bernd Simeon\footnote{TU Kaiserslautern, Dept. of Mathematics, Germany, simeon@mathematik.uni-kl.de},     
Utz Wever\footnote{Siemens AG, Technology, Munich , Germany, utz.wever@siemens.com}}

\begin{document}
\maketitle
\begin{abstract}
The paper presents a new reduction method designed for dynamic contact problems. Recently, we have proposed an efficient reduction scheme for the node-to-node formulation, that leads to Linear Complementarity Problems (LCP). Here, we enhance the underlying contact problem to a node-to-segment formulation. Due to the application of the dual approach, a Nonlinear Complementarity Problem (NCP) is obtained, where the node-to-segment condition is described by a quadratic inequality and is approximated by a sequence of  LCPs in each time step. These steps are performed in a reduced approximation space, while the contact treatment itself can be achieved by the Craig-Bampton method, which preserves the Lagrange multipliers and the nodal displacements at the contact zone. We think, that if the contact area is small compared to the overall structure, the reduction scheme performs very efficiently, since the contact shape is entirely recovered. The performance of the resulting reduction method is assessed on two 2D computational examples. 
\end{abstract}

\vspace{1cm}
{\bf keywords:} Dynamic contact, model order reduction, adjoint approach, nonlinear complementarity problem, digital twin technology
\vspace{1cm}

\section{Introduction}
\label{sec:intro}
During the last couple of years, "digital twins" appear to be one of the key concepts towards digitalization. They provide computer-aided assistance for real-world models allowing to monitor their state at specified positions or make prediction for the the future states \cite{benner2020model}. Model order reduction realizes digital twins by enabling real-time execution. It reduces the degrees of freedom  while keeping the essential physical properties of the model. Moreover it enables the reusability of such models for other scenarios.
An extensive literature is devoted to model order reduction applied for different kinds of partial differential equations, see e.g., \cite{ lohmann2004, lohmann2006,Benner2014,Farhat2014,Chaturantabut2009, bennerwillcox}. Model reduction methods are typically categorized between data-driven methods \cite{radermacher2016pod, peherstorfer2016data, swischuk2019projection} or physics-based methods such as balanced truncation \cite{benner2016frequency, besselink2014model, kurschner2018balanced}, Krylov subspace methods \cite{Bai2002,ilyas2016krylov, breiten2010krylov, benner2011recycling, wolf2013,bentbib2018computational} and modal reduction \cite{Benner2014}.

In this paper we focus on reduction approaches for dynamic contact problems in linear elasticity. Such problems are naturally nonlinear due to the unknown moving contact interface, see, e.g., \cite{Fran75,KiOd88,wriggers2004computational,firrone2019} for introduction to contact mechanics
and \cite{krenciszek2014model} for first results on reduction
approaches in this matter. In \cite{Charbel2017}, a dynamic problem with a linear contact condition is discussed and both the displacements and the Lagrange multipliers are reduced  using methods such as the singular value decomposition or the non-negative matrix compression. In contrast, a reduction scheme for linear node-to-node contact problems was introduced in \cite{manvelyan2021efficient}, which reduces only the primal displacements. Such contact conditions allow, however, only sticking at the contact and loosing the contact in the normal direction . 

Mechanical contact problems represent variational problems under unilateral constraints \cite{wriggers2004computational, capatina2014variational}, which are commonly solved by the penalty method \cite{hunvek1993penalty} or the augmented Lagrange method \cite{Gill81} after the discretization. Despite their advantages, these methods turn out to be inefficient for designing reduction schemes. Moreover, the contact shape is usually not known a priori. This is a major challenge for reduction schemes, since the reduced model needs to account for the contact impact. 

In this work we propose a novel reduction method for the class of node-to-segment contact problems as described in, e.g., \cite{wriggers2004computational}. This discretization technique of the contact condition is widespread in the field of contact mechanics since it can handle large deformations and allows more flexibility at the contact zone such as sliding, see \cite{zavarise2009modified, zavarise1998segment, zavarise2009node}.
Solvers for such contact problems have been integrated into nonlinear finite element commercial software \cite{nastran,ansys}. For early implementations of node-to-segment techniques we refer to \cite{hallquist1979nike2d,hughes1977finite}.

The novel reduction approach is purely physics-informed, i.e., no trajectories of the full model are required. The projection matrix can be computed based on the system matrices and the external load vector. Once the reduced model is established  within an offline procedure, we turn to the adjoint space of the reduced displacement, i.e., the space of the dual Lagrange multipliers. In contrast to \cite{manvelyan2021efficient}, the time-dependent dual problem represents a Nonlinear Complementarity Problem (NCP). The latter is approximated by a sequence of Linear Complementarity Problems (LCPs) which can be solved within fixed-point iterations for each time step. In terms of convergence, a quadratic convergence rate for the fixed-point iterations is guaranteed due the general theoretical 
framework \cite{adly2016newton,pang1982iterative,josephy1979newton}. Combining with the fact that the sequential LCPs stem from the reduced space and the number of the constraints is small, the convergence rate provides major savings with respect to execution time.
 
More specifically, the model order reduction is essentially based on a partitioning of the  the overall displacement degrees of freedom in
the flavor of the Craig-Bampton method \cite{craigbampton}, which separates the interior nodes, i.e., the slave nodes, from the contact nodes, i.e., the master nodes. Due to the nonlinear constraints, the system stiffness matrix turns out to depend on the Lagrange multipliers. However, the partitioning of the nodes provides linear system block-matrices for the slave nodes allowing an accurate reduction. After the reduced slave nodes are computed, the nonlinear block-matrices of the master nodes are utilized for the static condensation restoring the coupling between the master and slave nodes. Thus the partitioning provides essential improvement in the accuracy for recovering the Lagrange multipliers and the contact shape as well. To put it in a nutshell, our proposed reduction scheme combines  Krylov subspace methods, adjoint methods for the linearized problem and the Craig-Bampton partitioning technique. 
The overall methodology is general in nature, but  so far it has been worked out in detail and tested in the planar case only.

Our current contribution presents the new reduction algorithm in the following steps: In Section \ref{sec:models} the variational formulation of dynamic contact problems with quadratic constraints is derived. The common node-to-segment technique is described in Section \ref{sec:contact}. Furthermore, Section \ref{sec:staticncp} and Section \ref{sec:dynamicncp} discuss the solving procedure of the contact problem. Instead of solving the underlying NCP, the corresponding linearized system (LCP) is sequentially solved within fixed-point iterations. The resulting reduction method is outlined in Section \ref{sec:reduction}. The performance of the approach is highlighted by means of two numerical examples in Section \ref{sec:applications}. Finally, a conclusion is drawn in Section \ref{sec:conclusions}.

\section{The Generic Formulation} 
\label{sec:models}
In this section we first outline the governing equations for dynamic contact in the semi-discretized form. Afterwards, the derivation of the variational formulation is briefly sketched.

\subsection{The Generic Semi-discretized Model}
In this paper we restrict ourselves to finite element methods  for a comprehensive contact treatment. The contact is discretized by the node-to-segment technique that allows the sliding of a node over segments defined on the contact area.
As explained below, such an approach leads to a quadratic inequality in the displacement vector that is combined with the equations of structural dynamics. In this way, a generic model for the system of semi-discretized equations can be derived that reads
\begin{subequations}\label{eq:KKT_transient1}
\begin{eqnarray}
	\Mat{M} \Vec{\ddot{q}} + \Mat{K} \Vec{q} -  \big(  (\Mat{D}+\Mat{D}^T) : \Vec{q}\big)^T\Vec{\lambda} - \Vec{f}  - \Mat{C}^T \Vec{\lambda} & = & 0, \label{ode1} \\ 
	(\Mat{D}:\Vec{q})\Vec{q}+\Mat{C} \Vec{q} + \Vec{b} & \geq & 0, \label{const11} \\
	\Vec{\lambda} & \geq & 0 \label{const12},\\
	\Vec{\lambda}^T\big((\Mat{D}:\Vec{q})\Vec{q}+\Mat{C} \Vec{q} + \Vec{b}\big) & = & 0, \label{const13}
\end{eqnarray}
\end{subequations}
where $\Vec{q}(t) \in\R^N$ denotes the vector of nodal displacement variables, $\Mat{M}, \Mat{K} \in \R^{N\times N} $ the system mass and stiffness matrix, respectively, and $\Mat{f}(t) \in \R^N$ the external loads. 
Moreover, $\Mat{D} \in \R^{m\times N\times N} $ is a third-order tensor,  $\Mat{C}\in \R^{m\times N}$ stands for the constraint matrix and
 $\Vec{b} \in \R^m$ for a given initial clearance vector.
In \eqref{const11} and \eqref{const13} the quadratic term is defined as 
$[ (\Mat{D}:\Vec{q})\Vec{q}]_k
= \Vec{q}^T \Mat{D}_k \Vec{q},$ with $\Vec{D}_k \in \R^{N \times N}$ for $k=1,\dots,m$.
Analogously, the inequality constraints are to be interpreted component-wise and
 stand for the discretized non-penetration conditions of all variables in the contact interface where the vector of Lagrange multipliers $\Vec{\lambda}(t) \in \R^m$ 
 enforces the constraint in the dynamic equation \eqref{ode1}.
 The multiplier $\Vec{\lambda}$, standing for the discretized contact pressure, is always positive and its inner product with the constraints satisfies a complementarity condition.
 Note also that 
 \begin{equation}\label{eq:structconstr}
  ((\Mat{D}+\Mat{D}^T) : \Vec{q}\big)^T\Vec{\lambda} = \sum_{k=1}^m(\Mat{D}_k + \Mat{D}_k^T)\Vec{q}\lambda_k
  = \left( \sum_{k=1}^m\lambda_k(\Mat{D}_k + \Mat{D}_k^T) \right) \Vec{q},
 \end{equation}
 which means that this term can be viewed as a nonlinear contribution to the stiffness matrix that depends on the Lagrange multipliers.
 
In this paper, we will introduce a physics-based model order reduction method for the semi-discretized equations (\ref{eq:KKT_transient1}) that projects
the nodal variables $\Vec{q}$ onto a much smaller space while explicitly preserving
the constraints and the Lagrange multipliers. For this purpose, we assume a frictionless, adhesive-free normal contact in combination with small deformation theory and a linear-elastic material. Those assumptions will lead to 
(\ref{eq:KKT_transient1}) if 
\begin{itemize}
\item[(i)] the Lagrange multiplier method is used to enforce the non-penetration condition 
\item[(ii)] the finite element method is applied for the spatial discretization
\item[(iii)] the node-to-segment contact technique is used for the discretization of the contact zone.
\end{itemize}
We do not consider approaches like the Augmented Lagrangian method or the Nitsche method. The same holds for additional friction effects. In contrast to  \cite{manvelyan2021efficient}, the node-to-segment provides quadratic inequalities as constraints. This technique will be sketched in greater detail in Section \ref{sec:contact}.

For a contact problem with $k$ bodies that fits into the framework described so far,
the mass matrix $\Vec{M}$ and the stiffness matrix $\Vec{K}$ consist of $k$ diagonal block matrices that stem from the discretization of each individual body. The constraint tensors $\Mat{D}$ and $\Vec{C}$ will then typically exhibit a sparse structure where each row stands for a node-to-segment contact condition.
For further details of various contact models and discretization schemes we refer to \cite{KiOd88,wriggers2004computational}. Before we proceed further with the numerical approaches, we briefly summarize the setting of the classical variational formulation of  contact problems in the dynamic case.

\subsection{Strong Formulation of the Dynamic Contact Problem}

We introduce some notational preliminaries. The initial undeformed solid  $\Omega \subset \R^d, d=2$ or $d=3$ is a disjoint union of bounded connected subdomains $\Omega_i,\ i \in \mathcal{I},$ i.e., $\Omega = \bigcup_{i \in  I}\Omega_i$ with $\Omega_i \cap \Omega_j = \emptyset$ for $ i \neq j.$ On $\Omega$ we consider a frictionless, adhesive-free normal multi-body contact problem and a linear-elastic material describing the displacement field $\Vec{u}(\Vec{x},t)\in \mathbb{R}^d$ and the contact pressure 
 $p(\Vec{x},t)\in \R$ where $x \in \Omega$ is the spatial variable 
 and $t \in [t_0,T]$ stands for the temporal variable. The restriction of the displacement field on the corresponding domain we denote as $\Vec{u}_i = \Vec{u}_{|\Omega_i}.$
 The boundary of each elastic body is decomposed into
$\partial\Omega_i = \Gamma^i_D\cup\Gamma^i_N\cup\Gamma^i_C,$ where the latter represents the contact interface. The contact problem in strong form is then given by a dynamic boundary value problem (BVP):
For all $i \in \mathcal{I}$ and for the corresponding initial data $\Mat{u}_i(\cdot,t_0) = \Mat{u}^i_0,\ \Mat{u}_{i,t}(\cdot,t_0) = \Mat{v}^i_0 $ find $\Vec{u}_i : \Omega_i \to \mathbb{R}^d$ such that
\begin{equation}\label{dbvp}
\begin{aligned} 
\rho_i\Mat{u}_{i,tt}-\mbox{div} \Mat{\sigma}(\Vec{u}_i) &= \Vec{F}_i \quad \mbox{in } \Omega_i \times [t_0,T],\\
\Vec{u}_i &=0  \quad \mbox{on } \Gamma^i_{D} \times [t_0,T],\\
\Mat{\sigma}(\Vec{u}_i)\cdot\Vec{n}_i&=\Vec{\tau}_i \quad \mbox{on } \Gamma^i_{N} \times [t_0,T],
\end{aligned}
\end{equation}
where the constant $\rho_i$ denotes the mass density of the body $\Omega_i$, $\Vec{F}_i(\Vec{x},t) \in \R^d$ the volume force and
$\Vec{\tau}_i(\Vec{x},t) \in \R^d$ the surface traction. 
The stress tensor
$\Mat{\sigma}(\Vec{u}_i) \in \R^{d\times d}$ and the linearized
strain tensor $ \Mat{e}(\Vec{u}_i) \in \R^{d\times d}$ are  given by
\begin{align}\label{stress_strain}
  \Mat{\sigma}(\Vec{u}_i) &= \frac{E_i}{1+\nu_i} \Mat{e}(\Vec{u}_i) + \frac{\nu_i E_i}{(1+\nu_i)(1-2 \nu_i)} \mbox{trace}(\Mat{e}(\Vec{u}_i)) \Mat{I},\\
  \Mat{e}(\Vec{u}_i) &= \frac{1}{2}(\nabla \Vec{u}_i+\nabla \Vec{u}_i^T)
\end{align}
with Young's modulus $E_i \geq 0$ and Poisson's ratio $-1\leq \nu_i \leq 0.5.$
Furthermore, the boundary value problem \eqref{dbvp} is subject to the contact conditions 
\begin{equation}\label{bcDN} 
g_i \geq 0, \quad
p_i \geq 0, \quad g_i \cdot p_i = 0  \quad \mbox{on } \Gamma^i_C \times [t_0,T].
\end{equation}
with $g_i = g_{|\partial\Omega_i},$ where $g(\cdot, \Vec{u}(\cdot)):\partial \Omega \to \R $ is the gap function described as
\begin{equation}\label{gap_function}
g_{|\partial\Omega_i}(\Vec{x}_i, \Vec{u}_i) = \min_{j\in \mathcal{I}/\{i\}}\big(\Vec{\xi}_j(\Vec{x}, \Vec{u}) - (\Vec{x}_i + \Vec{u}_i(\Vec{x}_i,t)\big)^T \Vec{n}_j(\Vec{\xi}_j), \quad \Vec{x}_i \in \partial \Omega_i
\end{equation}
and $p_i$ is the contact pressure on $\partial\Omega_i$, respectively. Note, that $\Vec{n}_j$ is the outer normal vector on $\partial\Omega_j$ and $\Vec{\xi}_j(\Vec{x}_i, \Vec{u}_i)$ is the projection point of $\Vec{x}_i$ to the body $\Omega_j$ in the outer normal direction, fulfilling
\begin{align}\label{xsi}
    \Vec{\xi}_j(\Vec{x}_i, \Vec{u}_i) = \argmin_{\Vec{\xi} \in \partial\Omega_j(t)}\|(\Vec{x}_i + \Vec{u}_i(\Vec{x}_i,t))-\Vec{\xi}\|, \quad \Vec{x}_i \in \partial\Omega_i.
\end{align}
On the one hand the gap function \eqref{gap_function} simply measures  the distance between $\Vec{x}_i \in \partial\Omega_i$ and $\Vec{\xi}_j \in \partial\Omega_j.$
But on the other hand depending on the underlying geometry, \eqref{gap_function} can lead to quite complex expressions. Moreover, in case of the non-uniqueness of \eqref{xsi}, the gap function \eqref{gap_function} can turn into a non-smooth function. For simplicity, we assume that there exist functionals $\Mat{H}:\partial\Omega \to \R^{d \times d}$ and $\Mat{\eta}:\partial\Omega \to \R^d$ such that for all $i\in\mathcal{I}$ it holds
\begin{align}\label{normal}
    \Vec{n}_i = \Mat{H}(\Vec{x}_i)\Vec{u}_i + \Vec{\eta}(\Vec{x}_i), \quad \Vec{x}_i \in \partial\Omega_i,
\end{align}
i.e., the outer normal $\Vec{n}_i$ depends on the displacement vector $\Vec{u}_i$ affine-linearly.

In \cite{manvelyan2021efficient} we considered a non-penetration condition allowing only normal contact such that no movement in the tangential direction was possible, see also \cite{wriggers2004computational}. In particular the gap function $g$ was linear. In this work, however, using the assumption \eqref{normal} we want to generalize the constraint to a quadratic condition allowing also sliding movement such that a point can move along the contact area.
\begin{figure}[H]
\centering
\includegraphics[scale=0.5]{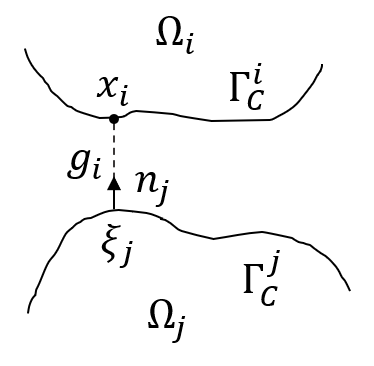}
\caption{Two-body contact problem}
\end{figure}
\subsection{Weak Formulation of the Dynamic Contact Problem}

In order to state the weak formulation of the contact problem, we introduce the function spaces
\begin{eqnarray}
\mathcal{V} &:= & \left\{ \Vec{v} \in H^1(\Omega)^d :
                           \Vec{v} = 0 \, \mbox{on } \Gamma_D \right\}, \label{defV} \\
\mathcal{L} &:= & \left\{ \mu \in H^{1/2}(\Gamma_C)^\prime : 
                         \int_{\Gamma_C} \mu w \, ds \geq 0 \,\, \forall
                          w \in H^{1/2}(\Gamma_C), w \geq 0 \right\},
                          \label{defL}
\end{eqnarray}
where $\mathcal{V}$ is the set of admissible displacements and is decomposed body-wise such that for all $v \in \mathcal{V}$ it holds $v_{|\Omega_i} \in H^1(\Omega_i)^d $ and  $v_{|\Omega_i} = 0 $ on $\Gamma^i_D$. The set $\mathcal{L}$ is the convex cone of the Lagrange multipliers defined on $\Gamma_C,$ which is the union of all $\Gamma^i_C.$  Furthermore, the following abstract notation is introduced:
\begin{equation}\label{defANot}
\langle \rho \Vec{u}_{tt}, \Vec{v} \rangle := \int_\Omega \Vec{v}^T \rho \Vec{u}_{tt} dx, \,\,
a(\Vec{u}, \Vec{v}) := \int_\Omega\Mat{\sigma}(\Vec{u}):\Mat{e}(\Vec{v})dx,
\,\,
\langle \ell , \Vec{v} \rangle := \int_\Omega \Vec{v}^T \Mat{F}dx + \int_{\Gamma_N} \Vec{v}^T\Vec{\tau}ds.
\end{equation}
 The non-penetration condition in weak form with test function $\mu \in {\cal L}$ is in general given by 
\begin{equation}\label{contactweak}
 0 \leq \int_{\Gamma_C} g \mu ds = 
        \int_{\Gamma_C} g (\mu-p) ds.
\end{equation}
In particular for a two-body problem, i.e. $\mathcal{I} = \{1,2\},$ the right-hand side of  \eqref{contactweak} reads:
\begin{equation}\label{twoweak}
\int_{\Gamma_C} g \mu ds = \int_{\Gamma_C}\big(\Vec{u}^T\Mat{H}\Vec{u} + \Vec{u}^T(\Vec{\eta} + \Mat{H}^T(\Vec{x}-\Vec{\xi})) + (\Vec{x}-\Vec{\xi})^T\Vec{\eta}\big)(\mu-p) ds.
\end{equation}

Note, that we provide the weak formulation solely for the  two-body problem. The extension to the general multi-body case is straightforward.  Using \eqref{twoweak} we define the trilinear form on ${\cal V} \times {\cal V} \times {\cal L}$
\begin{equation}\label{defTriFm}    
\mathit{d}(\Vec{v},\Vec{v},p) :=   \int_{\Gamma_C}  \Vec{v}^T\Mat{H}\Vec{v} \cdot p \,  ds,
\end{equation}
the bilinear form on ${\cal V} \times {\cal L}$
\begin{align}\label{defBilFm}    
\mathit{c}(\Vec{v},p) :=  \int_{\Gamma_C}  \Vec{v}^T(\Vec{\eta} + \Mat{H}^T(\Vec{x}-\Vec{\xi})) \cdot p \,  ds
\end{align}
and the linear form on $\cal L$ 
\begin{equation}\label{defLFm}
\langle \mathit{b}, p \rangle := \int_{\Gamma_C} (\Vec{x}-\Vec{\xi})^T\Vec{\eta} \cdot p \, ds.
\end{equation}
After these preparations, the weak form of the dynamic contact problem is stated as follows: For each $ t \in [t_0, T] $ find the displacement field $ \Vec{u}(\cdot,t) \in \mathcal{V} $ and the contact pressure $ p(\cdot,t) \in \mathcal{L} $ such that  
\begin{equation}\label{elast_weak}
\begin{aligned} 
&\langle \rho \Vec{u}_{tt}, \Vec{v} \rangle + a(\Vec{u}, \Vec{v}) - \mathit{d}(\Vec{u},\Vec{v},p)  = \langle \ell,\Mat{v} \rangle  + c(\Vec{v}, p) \quad &&\mbox{for all } \Vec{v} \in \mathcal{V},\\
& \mathit{d}(\Vec{u},\Vec{u}, \mu-p) + \mathit{c}(\Vec{u}, \mu-p) +  \langle \mathit{b}, \mu-p \rangle  \geq 0 \quad &&\mbox{for all }  \mu \in \mathcal{L}. 
\end{aligned}
\end{equation}
The system \eqref{elast_weak} is discretized with respect to the spatial variable by applying the standard Galerkin projection $\Vec{u}(\Vec{x},t) \doteq \sum_{i=1}^N \Vec{\phi}_i(\Vec{x}) q_i(t)
= \Mat{\Phi}(\Vec{x}) \Vec{q}(t)$ with basis functions $\Vec{\phi}_i$ and
nodal variables $\Vec{q} = (q_1, \ldots, q_N)^T \in \R^N$ to the displacement field.

Moreover, the integral over $\Gamma_C$ in the weak dynamic equation can be approximated by a discrete sum
\begin{equation}
\begin{aligned}\label{discretGammaC}
         \int_{\Gamma_C} \Vec{v}^T\Mat{H}\Vec{u} \cdot p \,  ds\ &+ \int_{\Gamma_C}  \Vec{v}^T(\Vec{\eta} + \Mat{H}^T(\Vec{x}-\Vec{\xi})) \cdot p \,  ds \\ &\doteq
                     \sum_{\ell=1}^{m} \Big(\varepsilon^1_\ell \Vec{v}(\Vec{x}_\ell)^T\Mat{H}\Vec{u}(\Vec{x}_\ell) + \varepsilon^2_\ell \Vec{v}(\Vec{x}_\ell)^T ( \Vec{\eta}(\Vec{x}_\ell)-\Mat{H}(\Vec{x}_\ell+\Vec{\xi}(\Vec{x}_\ell))\Big) \cdot p(\Vec{x}_\ell,t)
\end{aligned}
\end{equation}
with the constants $ \varepsilon^i_\ell = \varepsilon^i_\ell(A_\ell),\ i = 1,2 $ where $A_\ell$ is the area around node $\Vec{x}_{\ell}.$
Putting finally $\Vec{\lambda}(t) := (p(\Vec{x}_1,t), \ldots, p(\Vec{x}_{m},t))^T \in \R^m$ as discrete pressure variable,
the discretized constraint with the third-order tensor $\Mat{D}\in \R^{m \times N \times N},$  the matrix $\Mat{C} \in \R^{m \times N}$ and contact offset $\Vec{b}\in\R^m$ lead to a quadratic inequality of the form \eqref{eq:KKT_transient1}. Note that only a few displacements are involved for the contact condition and hence most of the entries of $\Mat{D}$ and $\Mat{C}$ are zero.

\section{The Node-To-Segment Contact Condition} 
\label{sec:contact}
In this section we study the contact condition using the node-to-segment technique in more detail. As already indicated above
it turns out that this condition is quadratic in terms of the displacements. In the following we assume a 2D setting where linear finite elements are used for the discretization. Let
\begin{equation}
	\textbf{p} = (\text{p}_1,\ \text{p}_2)^T \in \mathbb{R}^2, \quad 
	\tilde{\textbf{p}} = ( \tilde{\text{p}}_1,\ \tilde{\text{p}}_2)^T \in \mathbb{R}^2, \quad 
	\textbf{p} \neq \tilde{\textbf{p}}
\end{equation}
be two neighbouring contact nodes. They define a contact segment that an be written as
\begin{equation}
	\mathcal{S}(\textbf{p},\tilde{\textbf{p}}) = \left \{ \Vec{s}(t) = \textbf{p} + t(\tilde{\textbf{p}}-\textbf{p}) 
	\ \vert \ t \in [0,1] \right \}.
\end{equation}

\begin{figure}[hb]
	\begin{subfigure}{0.5 \textwidth}
		\centering
		\begin{tikzpicture}[scale = 0.4]
		\coordinate (p) at (-3,-2);
		\coordinate (q) at (6,4);
		\coordinate (r) at (5,-1);
		\coordinate (s) at (3,2);
		
		\draw[thick] (p) -- (q);
		\draw[densely dashed, thick,] (s) -- (r);
		
		\path (r) -- (s) -- (q) pic ["\Huge $\cdot$", draw, -] {angle=r--s--q};
		
		\draw [fill] (p) circle [radius=0.2] node [below=4] {$\textbf{p}$};
		\draw [fill] (q) circle [radius=0.2] node [below=4] {$\tilde{\textbf{p}}$};
		\draw [fill] (r) circle [radius=0.2] node [below=4] {$\Vec{r}$};
		\draw [fill] (s) circle [radius=0.2] node [left=4] {$\Vec{s}(t^*)$};
		
		\end{tikzpicture}
		\caption{$g_{\Vec{s}}(\Vec{r})>0$\label{fig:abstand} }
	\end{subfigure}%
	\begin{subfigure}{0.5 \textwidth}
		\centering
		\begin{tikzpicture}[scale = 0.4]
		\coordinate (p) at (-3,-2);
		\coordinate (q) at (6,4);
		\coordinate (r) at (5,-1);
		\coordinate (mr) at (3,-1);
		
		\draw[thick] (p) -- (q);
		\draw[ultra thick,->] (p) -- node[above] {} (r);
		\draw[ultra thick,->] (p) -- node[below,shift={(-2mm,0mm)}] {$\Vec{n}$} ($(p)!(mr)!90:(q)$);
		
		\path (p) -- ($(p)!(mr)!90:(q)$) coordinate (nd);
		
		\path (nd) -- (p) -- (r) pic [draw,thick,fill=gray,fill opacity=0.4,text opacity=1,angle radius=0.6cm,"{$\measuredangle (\Vec{n},\Vec{r}-\textbf{p})$}" shift={(12mm,-2mm)}] {angle=nd--p--r};

		\draw [fill] (p) circle [radius=0.2] node [left=4] {$\textbf{p}$};
		\draw [fill] (q) circle [radius=0.2] node [below=4] {$\tilde{\textbf{p}}$};
		\draw [fill] (r) circle [radius=0.2] node [above=4] {$\Vec{r}$};		
		\end{tikzpicture}
		\caption{\label{fig:measure>0} $\mathcal{M}(\textbf{p},\tilde{\textbf{p}},\Vec{r}) > 0$}
	\end{subfigure}

	\caption{Geometrical illustration}
\end{figure}
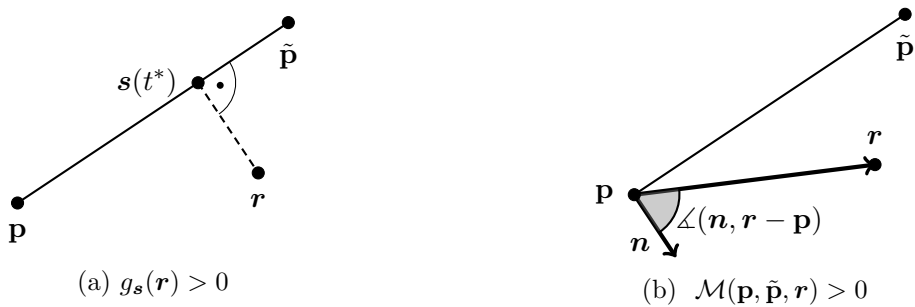

Furthermore, we assume that after the discretization a node-segment correspondance is established at the contact zone. Let $\Vec{r} \in \mathbb{R}^2$ be the corresponding node of the segment $\mathcal{S}$. The distance of the segment $\mathcal{S}$ to the node $\Vec{r}$ is given by (see Figure \ref{fig:abstand})
\begin{equation} \label{eq:distance}
	g_{\Vec{s}}(\Vec{r}) =\min_{t\in[0,1]} \| \Vec{s}(t) - \Vec{r} \|_2.
\end{equation}

Similarly as for the gap function in \eqref{gap_function}, the distance function in \eqref{eq:distance} can lead to a quite complex expression. Therefore we define an angular function representing an easier alternative to \eqref{eq:distance}. 
For a node $\Vec{r}$ let $\mathcal{S}(\textbf{p},\tilde{\textbf{p}})$ be the nearest segment to $\Vec{r}$ and let $\Vec{n}$ denote the normal vector orthogonal to the line containing $\mathcal{S}$. Then the angular function is given by (see Figure \ref{fig:measure>0})
\begin{equation} \label{eq:measure}
	\mathcal{M}(\textbf{p},\tilde{\textbf{p}},\Vec{r}) = \langle \Vec{n}, \Vec{r}-\textbf{p} \rangle = \Vec{n}^T(\Vec{r}-\textbf{p})
\end{equation}
Note, that for two nonzero vectors $\Vec{v}_1, \Vec{v}_2 \in \R^\ell \setminus \{ 0 \}$ the equality $	\langle \Vec{v}_1, \Vec{v}_2 \rangle = \| \Vec{v}_1 \|_2 \cdot \| \Vec{v}_2 \|_2 \cdot \cos(\measuredangle (\Vec{v},\Vec{w}))$ always holds, where $\measuredangle (\Vec{v},\Vec{w}) \in [0,\pi]$ denotes the angle enclosed by $\Vec{v}$ and $\Vec{w}.$ 

For \eqref{eq:measure} three different scenarios can be distinguished:
\begin{itemize}
	\item $\mathcal{M}(\textbf{p},\tilde{\textbf{p}},\Vec{r}) > 0$, which means that $\Vec{r}$ lies in the half of the line with positive normal.
	\item $\mathcal{M}(\textbf{p},\tilde{\textbf{p}},\Vec{r}) = 0$, which means that $\Vec{r}$ lies on the line containing the segment $\mathcal{S}$.
	\item $\mathcal{M}(\textbf{p},\tilde{\textbf{p}},\Vec{r}) < 0$, which means that $\Vec{r}$ lies in the half of the line with negative normal.
\end{itemize}

The normal $\Vec{n}$ to $\mathcal{S}$ depends on the end nodes of $\mathcal{S}$ and is given by
\begin{equation}
	\Vec{n}(\textbf{p}, \tilde{\textbf{p}}) = 
	\begin{pmatrix}
		0 & -1 \\
		1 & 0
	\end{pmatrix}
	\left ( \tilde{\textbf{p}} - \textbf{p} \right ).
\end{equation}

In order to derive a formulation for \eqref{eq:measure} in terms of the displacements and also to keep the notation concise, we write $\Vec{u}^i_0 + \Vec{u}^i$ with the vectors $\Vec{u^i}_0 \in \mathbb{R}^2$ standing for the initial position of the node $i$ and $\Vec{u}^i \in \mathbb{R}^2 $ for the corresponding displacement. The label $i$ runs through $\{ \textbf{p},\tilde{\textbf{p}},\Vec{r} \}.$ Based on this notation, the function \eqref{eq:measure} can be reformulated as 

\begin{equation} \label{ContactCondition}
\begin{array}{ccrcl}
\mathcal{M}(\textbf{p},\tilde{\textbf{p}},\Vec{r}) & = & 
	\left( \Vec{u}^{\tilde{\textbf{p}}}_0 - \Vec{u}^{\textbf{p}}_0 \right)^T & 
	\begin{pmatrix} 0 & -1 \\ 1 & 0 \end{pmatrix} & 
	\left( \Vec{u}^{\Vec{r}}_0 - \Vec{u}^{\textbf{p}}_0 \right) \\
\ & + & 
	\left( \Vec{u}^{\tilde{\textbf{p}}}_0 - \Vec{u}^{\textbf{p}}_0 \right)^T & 
	\begin{pmatrix} 0 & -1 \\ 1 & 0 \end{pmatrix} & 
	\left( \Vec{u}^{\Vec{r}} - \Vec{u}^{\textbf{p}} \right) \\
\ & + & 
	\left( \Vec{u}^{\Vec{r}}_0 - \Vec{u}^{\textbf{p}}_0 \right)^T & 
	\begin{pmatrix} 0 & 1 \\ -1 & 0 \end{pmatrix} & 
	\left( \Vec{u}^{\tilde{\textbf{p}}} - \Vec{u}^{\textbf{p}} \right) \\
\ & + & 
	\left( \Vec{u}^{\tilde{\textbf{p}}} - \Vec{u}^{\textbf{p}} \right)^T & 
	\begin{pmatrix} 0 & -1 \\ 1 & 0 \end{pmatrix} & 
	\left( \Vec{u}^{\Vec{r}} - \Vec{u}^{\textbf{p}} \right).
\end{array}
\end{equation}

The function $\mathcal{M}$ defined in \eqref{ContactCondition} is at most of quadratic order in terms of the displacements. 
Furthermore, by assembling the corresponding terms of the summation within the equation \eqref{ContactCondition} in a proper way we can determine the global constraint matrices $\Mat{D}_k \in \R^{n \times n}$, the global vectors $\Vec{c}_k \in \R^n$ and and the scalars $b_k \in \R$. For the nodal displacement vector $\Vec{q} \in \mathbb{R}^n$ and $m$ contact nodes the contact condition reads
\begin{equation} \label{ContactConst.}
	\Vec{q}^T \Mat{D}_k  \Vec{q} + \Vec{c}_k^T \Vec{q} + b_k \geq 0, \quad k=1, \dots, m
\end{equation}
In a more compact form (\ref{ContactConst.}) can be stated as 
\begin{equation} \label{QuadraticContactConst.}
	\begin{pmatrix}
	\Vec{q}^T \Mat{D}_1 \Vec{q} \\
	\vdots \\
	\Vec{q}^T \Mat{D}_m \Vec{q}
	\end{pmatrix}
	+ \Mat{C} \Vec{q} + \Vec{b} \geq 0
\end{equation} 
with 
\begin{equation}
	\Mat{C} = 
	\begin{pmatrix}
	\Vec{c}_1^T \\
	\vdots \\
	\Vec{c}_m^T
	\end{pmatrix} \in \R^{m \times n} \quad \text{and} \quad 
	\Vec{b} = 
	\begin{pmatrix}
	b_1 \\
	\vdots \\
	b_m 
	\end{pmatrix} \in \R^m
\end{equation}
The latter represents a quadratic contact condition derived from the node-to-segment discretization technique.

\section{Solving Static Contact Problems} 
\label{sec:staticncp}
A common approach to solve mechanical contact problems is the augmented Lagrangian method \cite{Powell69,Gill81,wriggers2004computational}. The method is very flexible to different descriptions of the contact region.
However, here we apply a dual approach, allowing us to turn to the adjoint problem resulting into a NCP. In the context of model order reduction for constrained problems dual approaches have a great benefit, which will be highlighted in the following sections. First we address the dual formulation and the solving procedure of a NCP in case of static contact problems.

\subsection{NCP for the static contact problem}
In the following we consider the quadratic contact condition
\begin{equation} 
\begin{pmatrix}
\Vec{q}^T \Mat{D}_1 \Vec{q} \\
\vdots \\
\Vec{q}^T \Mat{D}_m \Vec{q}
\end{pmatrix}
+ \Mat{C} \Vec{q} + \Vec{b} \geq 0.
\end{equation} 

Then, the variational formulation of the mechanical contact problem starts by considering the energy expression
\begin{equation}\label{eq:var-quad}
\min_{\Vec{q}\in\R^n} \frac{1}{2}\Vec{q}^T \Mat{K} \Vec{q} - 
\Vec{q}^T \Vec{f} - \sum_{k=1}^{m} \lambda_{k} \left( \Vec{q}^T \Mat{D}_k  \Vec{q} + \Vec{c}_k^T \Vec{q} + b_k \right).
\end{equation}
The corresponding KKT-conditions read
\begin{equation}\label{eq:KKT}
\begin{aligned} 
\Mat{K} \Vec{q} - \Vec{f} - \left( \sum_{k=1}^{m} \lambda_{k} (\Vec{D}_k + \Vec{D}_k^T) \right) \Vec{q} - \Mat{C}^T \Vec{\lambda} &= 0, \\
\Vec{q}^T \Mat{D}_k \Vec{q} + \Vec{c}_k^T \Vec{q} + b_k &\geq 0,\quad k=1,\dots,m, \\ 
\Vec{\lambda} &\geq 0, \\
\sum_{k=1}^{m} \lambda_{k} \left( \Vec{q}^T \Mat{D}_k \Vec{q} + \Vec{c}_k^T \Vec{q} + b_k \right) &= 0.
\end{aligned}
\end{equation}
Solving for the displacements $\Vec{q}$ in the KKT-conditions leads to the representation
\begin{equation} \label{eq:disp}
\Vec{q} = \left( \Mat{K} - \sum_{k=1}^{m} \lambda_{k} (\Mat{D}_k + \Mat{D}_k^T) \right)^{-1} 
(\Vec{f} + \Mat{C}^T \Vec{\lambda}).
\end{equation}

Next, we introduce the abbreviations
\begin{equation}\label{eq:abbK}
	\Mat{S}(\Vec{\lambda}) = \Mat{K} - \sum_{k=1}^{m} \lambda_{k} (\Mat{D}_k + \Mat{D}_k^T)
\end{equation}
and
\begin{equation} \label{static_w}
	\Vec{w}(\Vec{\lambda}) = \Mat{S}(\Vec{\lambda})^{-1} (\Vec{f} + \Mat{C}^T \Vec{\lambda})
\end{equation}
and insert (\ref{eq:disp}) into the second part of the KKT-conditions. This yields
\begin{equation} \label{eq:KKTL}
	\begin{array}{rcl}
	\Vec{w}(\Vec{\lambda})^T \Mat{D}_k \Vec{w}(\Vec{\lambda}) + \Vec{c}_k^T \Vec{w}(\Vec{\lambda}) + b_k & \geq & 0, \quad k = 1, \dots, m, \\[1mm]
	\Vec{\lambda} & \geq & 0,  \\
	\displaystyle
	\sum_{k=1}^{m} \lambda_{k} \left( \Vec{w}(\Vec{\lambda})^T \Mat{D}_k \Vec{w}(\Vec{\lambda}) + \Vec{c}_k^T \Vec{w}(\Vec{\lambda}) + b_k \right) & = & 0,
	\end{array}
\end{equation}
which is a nonlinear complementarity problem (NCP) for $\Vec{\lambda}$. After computing the Lagrange multiplier $\Vec{\lambda}$, the displacements are given by equation (\ref{eq:disp}).

\subsection{Solving a sequence of LCPs for the static contact problem}
A well-known and successful technique of solving NCPs is to consider their linear approximation, representing an LCP system, for which an LCP solver can be applied repeatedly. This method has been introduced in \cite{josephy1979newton} where the problem statement is given in a generalized framework. More about iterative techniques for solving variational inequalities can be found in \cite{pang1982iterative,adly2016newton}. A general form of an LCP is defined as follows:
\begin{equation} \label{eq:solveLCP1}
	\begin{array}{rcc}
	\Vec{B} + \Mat{A} \Vec{\lambda} & \geq & 0, \\
	\Vec{\lambda} & \geq & 0, \\
	\Vec{\lambda}^T (\Vec{B} + \Mat{A} \Vec{\lambda}) & = & 0.
	\end{array}
\end{equation}
LCP problems may be solved by applying either Fischer-Burmeister type algorithms \cite{Fischer92} or
more efficiently Lemke's algorithm (see \cite{LinearComplementarityProblem}, 4.4.5).  A  detailed discussion of Lemke's algorithm as well as a much more comprehensive study of Linear Complementarity Problems and can be found in \cite{LinearComplementarityProblem}. Specifically, we use a Python implementation of the algorithm \cite{lemke}.\\
The core idea, which below is also used for our reduction technique, is to exchange the original NCP problem
\begin{equation} \label{NCP_static}
	\Vec{\lambda} \geq \Vec{0}, \quad \Vec{F} (\Vec{\lambda}) \geq \Vec{0}, \quad \Vec{\lambda}^T \Vec{F} (\Vec{\lambda}) = 0
\end{equation}
by the linear approximation 
\begin{equation} \label{LCP_static}
	\Vec{\lambda} \geq \Vec{0}, \quad \Vec{F} (\Vec{\lambda}_0) + \nabla \Vec{F} (\Vec{\lambda}_0) (\Vec{\lambda} - \Vec{\lambda}_0) \geq \Vec{0}, \quad \Vec{\lambda}^T \left(\Vec{F} (\Vec{\lambda}_0) + \nabla \Vec{F} (\Vec{\lambda}_0) (\Vec{\lambda} - \Vec{\lambda}_0) \right) = 0,
\end{equation}
where 
\begin{equation}
	\Vec{F}: \R^m \to \R^m, \quad \Vec{F}(\Vec{\lambda}) = \big(\Vec{F}_1(\Vec{\lambda}), \ldots,\Vec{F}_m(\Vec{\lambda})\big)^T
\end{equation}
with
\begin{equation}
	\Vec{F}_k(\Vec{\lambda}) = \Vec{w}(\Vec{\lambda})^T \Mat{D}_k \Vec{w}(\Vec{\lambda}) + \Vec{c}_k^T \Vec{w}(\Vec{\lambda}) + b_k, \quad k = 1,\ldots,m
\end{equation}
represents the quadratic constraint function. The system \eqref{LCP_static} is solved by an LCP solver (e.g., the Lemke method) until a carefully chosen stopping criterion depending on the initial starting point $\Vec{\lambda}_0$ is satisfied. Then $\Vec{\lambda}_0$ is updated by $\Vec{\lambda}$ and
the LCP problem is solved again. Finally, this procedure yields an approximation for the solution vector $\Vec{\lambda}$ of the nonlinear problem \eqref{NCP_static}.

\section{Solving Dynamic Contact Problems} 
\label{sec:dynamicncp}
After the static contact problems have been addressed, we turn our attention to dynamic contact problems. Similar to the static case, the key idea for solving the time-dependent contact problems is to replace the nonlinear dual system by its linear approximation.

\subsection{Notational Convention}\label{notation}
For subsequent use, we introduce a new notation for our solution variables. By $\Vec{q}_{t_i}$ we denote $\Vec{q}(t_i),$ where $t_i$ stands for the $i-$th step of the discretized time grid. The Lagrange multiplier $\Vec{\lambda}^l_{t_i}$ carries two indices, where the superscript $l$ denotes the corresponding fixed-point iteration of the LCP. If it is clear from the context that the time instance $t_i$ is kept fixed, we make use of the  variable $\Vec{z}_l = \Vec{\lambda^l_{t_i}}.$

\subsection{NCP for the dynamic contact problem}
The dynamical model combined with the Kuhn-Tacker conditions \eqref{eq:KKT} leads to the system
\begin{subequations} \label{eq:KKT_transient}
\begin{eqnarray}
	&&\Mat{M} \Vec{\ddot{q}} + \Mat{K} \Vec{q} - \Vec{f} - \left( \sum_{k=1}^{m} \lambda_{k} (\Vec{D}_k + \Vec{D}_k^T) \right) \Vec{q} - \Mat{C}^T \Vec{\lambda} = 0, \label{ode} \\ 
	&&\Vec{q}^T \Mat{D}_k \Vec{q} + \Vec{c}_k^T \Vec{q} + b_k \geq 0,\quad k=1,\dots,m,\label{const1} \\[1mm]
	&&\Vec{\lambda} \geq 0, \label{const2} \\
	&&\sum_{k=1}^{m} \lambda_{k} \left( \Vec{q}^T \Mat{D}_k \Vec{q} + \Vec{c}_k^T \Vec{q} + b_k \right) = 0.\label{const3}
\end{eqnarray}
\end{subequations}
All four conditions in \eqref{eq:KKT_transient} must be fulfilled in each time step. The balance of momentum \eqref{ode} is discretized in time by the implicit Euler scheme, a straightforward and unconditionally stable method, also known as the first order backward differentiation formula. 
To this end, we replace the second order derivative by the finite difference
\begin{equation}\label{eq:findiff}
\ddot{\Vec{q}}(t+h) \approx \frac{1}{h^2}\left(\Vec{q}(t+h) - 2\Vec{q}(t) + \Vec{q}(t-h)\right),
\end{equation}
where $h$ is the time stepsize. Note that for the velocity it holds $\dot{\Vec{q}}(t+h) \approx
(\Vec{q}(t+h) - \Vec{q}(t))/h,$ a fact that is hidden behind \eqref{eq:findiff}. Furthermore, due to the presence of the
the second order time derivative, the implicit Euler leads to a two-step method.
Since \eqref{eq:findiff} possesses first order of accuracy only, 
 other integration schemes can be applied instead, e.g. the generalized-$\alpha $ method, which is very common in structural mechanics, see ~\cite{simeon}. In this work, however, we mainly focus on designing a reduction scheme and therefore, we continue working with the implicit Euler method.
Inserting \eqref{eq:findiff} into equation \eqref{eq:KKT_transient} and keeping in mind the abbreviation \eqref{eq:abbK}, we obtain
\begin{equation}\label{eq:system_disc}
\Mat{M}\left(\Vec{q}(t+h) -2\Vec{q}(t) + \Vec{q}(t-h)\right) + h^2\Mat{S}(\Vec{\lambda})\Vec{q}(t+h) = h^2\Vec{f}(t+h) + h^2\Mat{C}^T \Vec{\lambda}.  
\end{equation}
Assuming that the previous two time steps are known, \eqref{eq:system_disc} is solved for $\Vec{q}(t+h)$, which leads to
\begin{equation}\label{eq:system_disc_resolved}
\Vec{q}(t+h) = \left(\Mat{M}+h^2\Mat{S}(\Vec{\lambda})\right)^{-1}(h^2\Vec{f}(t+h) + h^2\Mat{C}^T\Vec{\lambda} +2\Mat{M}\Vec{q}(t) - \Mat{M}\Vec{q}(t-h)).
\end{equation}
We remark that the two-step time discretization requires initial values
$\Vec{q}(t_0)$ and $\Vec{q}(t_0+h)$ to start. If $\Vec{q}_0$ and
$\dot{\Vec{q}}_0$ as initial displacement and velocity are given, one can
compute $\Vec{q}(t_0+h) = \Vec{q}_0 + h \dot{\Vec{q}}_0$ by an explicit Euler step and then continue with the two-step formula \eqref{eq:system_disc}. The initialization of the Lagrange multiplier $\Vec{\lambda}$ is shortly addressed in Section \ref{sec:reduction}.
If the inequality constraints were replaced by the
equality constraints $\Vec{u}^T \Mat{D}_k \Vec{u} + \Vec{c}_k^T \Vec{u} + b_k = 0,\ k=1,\dots,m$, the index of the resulting differential-algebraic equation would equal 3, which means that special care 
must be taken for the time integration \cite{BrCP96,HaLR89,simeon}. Inserting \eqref{eq:system_disc_resolved} into the constraints \eqref{const1}-\eqref{const3} we obtain an NCP problem similar to \eqref{eq:KKTL}. After computing the Lagrange multiplier $\Vec{\lambda},$ the displacements are given by \eqref{eq:system_disc_resolved}. Note that the NCP has to be solved in each time step.

\subsection{Solving a sequence of LCPs for the dynamic contact problem}
After having introduced the transient model, we address the corresponding solution procedure in more detail.
The core idea  of the algorithm is the generalization of the iterative solver for the static problems to an iterative solver for the dynamic problemsby substituting the displacement vector in \eqref{static_w} by
\begin{equation}\label{eq:w}
\Vec{w}_{t_{i+1}}(\Vec{\lambda}_{t_{i+1}}) = \left(\Mat{M}+h^2\Mat{S}(\Vec{\lambda}_{t_i})\right)^{-1}(h^2\Vec{f}_{t_{i+1}} + h^2\Mat{C}^T\Vec{\lambda}_{t_{i+1}} +2\Mat{M}\Vec{w}_{t_{i}} - \Mat{M}\Vec{w}_{t_{i-1}}).
\end{equation}
The quadratic constraint function $\Vec{F}: \R^m \to \R^m$ is defined as 
\begin{align}\label{F_temp} 
	\Vec{F}_k(\Vec{\lambda}_{t_{i+1}}) &=\Vec{w}_{t_{i+1}}(\Vec{\lambda}_{t_{i+1}})^T \Mat{D}_k \Vec{w}_{t_{i+1}}(\Vec{\lambda}_{t_{i+1}}) + \Vec{c}_k^T \Vec{w}_{t_{i+1}}(\Vec{\lambda}_{t_{i+1}}) + b_k, \quad k = 1,\ldots,m.
\end{align}
The function \eqref{F_temp} defines for each time step the adjoint problem of \eqref{eq:KKT_transient}, which represents an NCP problem, i.e.,
\begin{equation} \label{NCP_dynamic}
	\Vec{\lambda}_{t_{i+1}} \geq \Vec{0}, \quad \Vec{F} (\Vec{\lambda}_{t_{i+1}}) \geq \Vec{0}, \quad \Vec{\lambda}_{t_{i+1}}^T \Vec{F} (\Vec{\lambda}_{t_{i+1}}) = 0.
\end{equation}
Instead of directly solving \eqref{NCP_dynamic}, we prefer to solve the linear approximation of \eqref{NCP_dynamic} within fixed-point iterations for each time step. For this purpose, we approximate the nonlinear constraint function \eqref{F_temp} by its linearization. As soon as the linearized constraint system comes into play, we make use of the abbreviation $\Vec{z}_l,$ see also Section \ref{notation}. This leads to the linearized system
\begin{equation}\label{NCP_linear_dynamic}
\begin{aligned} 
	\Vec{z}_l\geq \Vec{0}, \quad \Vec{F} (\Vec{z}_{l-1}) + \nabla \Vec{F}(\Vec{z}_{l-1}) (\Vec{z}_l - \Vec{z}_{l-1}) \geq \Vec{0},\\
	{\Vec{z}_l}^T \left(\Vec{F} (\Vec{z}_{l-1}) + \nabla \Vec{F} (\Vec{z}_{l-1}) (\Vec{z}_l - \Vec{z}_{l-1}) \right) = 0.
\end{aligned}
\end{equation}
Note that at the time step $t_{i+1}$ it holds for $\Vec{z}_0 = \Vec{\lambda}_{t_i},$ which is assumed to be known. For each $l>0$  the linearization is performed at $\Vec{z}_{l-1}.$ Furthermore, we write $\Mat{A}^{l-1} =\nabla \Vec{F}(\Vec{z}_{l-1})$ and $\Mat{B}^{l-1} =\Vec{F} (\Vec{z}_{l-1}) - \nabla \Vec{F}(\Vec{z}_{l-1})\Vec{z}_{l-1} $ such that the system \eqref{NCP_linear_dynamic} can be transformed into the LCP problem
\begin{equation} \label{eq:solveLCP2}
	\begin{array}{rcc}
	\Vec{B}^{l-1} + \Mat{A}^{l-1} \Vec{z}_{l} & \geq & 0, \\
	\Vec{z}_{l} & \geq & 0, \\
	\Vec{z}_{l}^T (\Vec{B}^{l-1} + \Mat{A}^{l-1} \Vec{z}_{l}) & = & 0.
	\end{array}
\end{equation}
The system \eqref{eq:solveLCP2} is repeatedly solved within a fixed-point iteration, until the error between $\Vec{z}_{l} $ and $\Vec{z}_{l-1}$ is small enough and the system converges, yielding the final $\Vec{z}_{l^*}$ for a certain $l^* >0.$\\ 
According to the general theoretical framework \cite{adly2016newton,pang1982iterative,josephy1979newton}, quadratic convergence is guaranteed when iteratively solving \eqref{NCP_linear_dynamic}. This means that only a few iterations suffice to obtain a solution for \eqref{NCP_dynamic} at each time step.  Within this work we do not pursue the theoretical aspects any further, but rather show a very good agreement with the theory based on our computational examples, see Section \ref{sec:applications}.

\subsection{Computation of the Jacobian}
In order to proceed with \eqref{NCP_linear_dynamic} we need to compute the Jacobian matrix of the function
\begin{equation}
	\Vec{F}: \R^m \to \R^m, \quad \Vec{F}(\Vec{\lambda}) = 
	\Vec{F} = \big(\Vec{F}_1(\Vec{\lambda}), \ldots,\Vec{F}_m(\Vec{\lambda})\big)^T
\end{equation}
with
\begin{equation}
	\Vec{F}_k(\Vec{\lambda}) = \Vec{w}(\Vec{\lambda})^T \Mat{D}_k \Vec{w}(\Vec{\lambda}) + \Vec{c}_k^T \Vec{w}(\Vec{\lambda}) + b_k.
\end{equation}
The derivative matrix is given by
\begin{equation}
	\Mat{DF}(\lambda) = 
	\begin{pmatrix}
	\frac{\partial}{\partial \Vec{\lambda}_1}\Vec{F}_1(\Vec{\lambda}) & \cdots & \frac{\partial}{\partial \Vec{\lambda}_m}\Vec{F}_1(\Vec{\lambda}) \\
	\vdots & \ddots & \vdots \\
	\frac{\partial}{\partial \Vec{\lambda}_1}\Vec{F}_m(\Vec{\lambda}) & \cdots & \frac{\partial}{\partial \Vec{\lambda}_m}\Vec{F}_m(\Vec{\lambda})
	\end{pmatrix} \in \R^{m \times m}.
\end{equation}
Due to the chain rule for the partial derivative of $\Vec{F}_i$ with respect to $\Vec{\lambda}_j$ it holds:
\begin{equation}\label{derivF}
\frac{\partial}{\partial \Vec{\lambda}_j}\Vec{F}_i(\Vec{\lambda}) = \left( \Vec{w}(\Vec{\lambda})^T(\Mat{D}_i + \Mat{D}_i^T) + \Vec{c}_i^T \right)\frac{\partial}{\partial \Vec{\lambda}_j}\Vec{w}(\Vec{\lambda}).
\end{equation}
In turn, the partial derivative of $\Vec{w}$ with respect to $\Vec{\lambda}_j$ results in
\begin{equation}\label{derivw}
\frac{\partial}{\partial \Vec{\lambda}_j}\Vec{w}(\Vec{\lambda}) = h^2(\Mat{M} +h^2\Mat{S}(\Vec{\lambda}))^{-1} \left( \Vec{c}_j + (\Mat{D}_j + \Mat{D}_j^T) \Vec{w}(\Vec{\lambda}) \right).
\end{equation}
Finally, inserting \eqref{derivw} into \eqref{derivF}, for $i,j \in \{ 1, \dots, m\}$ we obtain
\begin{equation}
\frac{\partial}{\partial \Vec{\lambda}_j}\Vec{F}_i(\Vec{\lambda}) = 
	h^2\left( \Vec{w}(\Vec{\lambda})^T(\Mat{D}_j + \Mat{D}_j^T) + \Vec{c}_j^T \right) (\Mat{M} +h^2\Mat{S}(\Vec{\lambda}))^{-T} \left( \Vec{c}_i + (\Mat{D}_i + \Mat{D}_i^T) \Vec{w}(\Vec{\lambda}) \right).
\end{equation}

For the sake of a more compact representation we introduce the term
\begin{equation}\label{abbrev}
	\Vec{Z}_k = \Vec{Z}_k(\Vec{\lambda}) = \Vec{w}(\Vec{\lambda})^T (\Mat{D}_k + \Mat{D}_k^T) + \Vec{c}_k^T, \quad k=1, \dots, m
\end{equation}
leading to the matrix
\begin{equation}
	\Mat{Z} = 
	\begin{pmatrix}
	\Vec{Z}_1 \\
	\vdots \\
	\Vec{Z}_m
	\end{pmatrix} \in \R^{m \times m}.
\end{equation}
Using the new abbreviation the partial derivatives can be computed as follows
\begin{equation}
	\frac{\partial}{\partial \Vec{\lambda}_j}\Vec{F}_i(\Vec{\lambda}) = h^2\Vec{Z}_i(\Mat{M}+ h^2\Mat{S}(\Vec{\lambda}))^{-1} \Vec{Z}_j^T.
\end{equation}
Finally, the Jacobian matrix can be stated as
\begin{equation}
	\Mat{DF}(\Vec{\lambda}) = h^2\Mat{Z}(\Mat{M}+h^2\Mat{S}(\Vec{\lambda}))^{-1} \Mat{Z}^T.
\end{equation}

\section{Reduction of the Contact Problem} 
\label{sec:reduction}
In many real world applications the numerical integration of constrained systems such as \eqref{eq:KKT_transient} is not possible in real time due to the large number $N$ of degrees of freedom. Model order reduction strategies \cite{Baur2014} introduce a reduced state $\Vec{q}_r\in\R^n$ with $n\ll N$, where $\Vec{q}_r$ is defined by
\begin{equation}\label{eq:reduction}
\Vec{q} = \Mat{Q}\Vec{q}_r,\quad \Mat{Q}\in\R^{N\times n}.
\end{equation}
\subsection{Computation of ROM}
One way to obtain the reduction matrix $\Mat{Q}$ for the system \eqref{ode} is to use modal reduction. Setting up the eigenvalue problem of \eqref{ode}
\begin{equation}
 \omega^2 \Mat{M} \Vec{v}=\Mat{K} \Vec{v} \label{eqn:ev problem}
\end{equation} 
and taking the first $n$ eigenvectors, the matrix $\Mat{Q}$ may be defined by
\begin{equation} \label{eq:Beigen}
\Mat{Q}=\left\{\Vec{v}_1,\cdots,\Vec{v}_n \right\}.
\end{equation} 
A preferable technique may be the Krylov subspace methods \cite{Bai2002,lohmann2006}. The subspace is defined by
\begin{equation}\label{eq:Bkrylov}
\Mat{Q}=\left\{\Mat{K}^{-1}\Vec{f}, 
\Mat{K}^{-1}\Mat{M} \Mat{K}^{-1}\Vec{f},\cdots,
(\Mat{K}^{-1}\Mat{M})^{n-1} \Mat{K}^{-1}\Vec{f} \right\}.
\end{equation}
The Krylov base may be computed by the Arnoldi algorithm, which delivers an orthonormal base of the subset.
Inserting the reduction \eqref{eq:reduction}, where $\Mat{Q}$ is defined by \eqref{eq:Beigen} or \eqref{eq:Bkrylov} into the differential equations system \eqref{ode} and multiplying by $\Mat{Q}^T$, one obtains the \emph{reduced  system}
\begin{equation}\label{eq:reduced_equation}
\hat{\Mat{M}} \ddot{\Vec{q}}_r + \hat{\Mat{S}}(\Vec{\lambda})\Vec{q}_r =
\Mat{Q}^T\left(\Vec{f}(t)+\Mat{C}^T \Vec{\lambda}(t)\right)
\end{equation}
where 
\begin{align}
\hat{\Mat{M}}&=\Mat{Q}^T \Mat{M}\Mat{Q} \in \R^{n\times n}, \\
\hat{\Mat{S}}(\Vec{\lambda}) &= \Mat{Q}^T\Mat{S}(\Vec{\lambda})\Mat{Q} = \Mat{Q}^T\Mat{K}\Mat{Q} - \sum_{k=1}^{m} \lambda_{k} \Mat{Q}^T(\Mat{D}_k + \Mat{D}_k^T)\Mat{Q} \in \R^{n\times n}
\end{align}
are the reduced mass and stiffness matrices.  
Similar to the full dimensional case, we apply implicit Euler to the reduced system \eqref{eq:reduced_equation}, which leads to
\begin{equation}\label{eq:system_time_red}
\Vec{q}_r(t+h) = \left(\hat{\Mat{M}}+h^2 \hat{\Mat{S}}(\Vec{\lambda})\right)^{-1}
\left(h^2 \Mat{Q}^T\Vec{f}(t+h) + h^2\Mat{Q}^T\Mat{C}^T\Vec{\lambda} +
 2\hat{\Mat{M}}\Vec{q}_r(t) - \hat{\Mat{M}}\Vec{q}_r(t-h))\right).
\end{equation}
Note that \eqref{eq:system_time_red} involves only operations with matrices and vectors in smaller dimensions $n$ or $m$.
Although $\Mat{Q}\in\R^{N\times n}$, the acting force $\Vec{f}(t+h)$  has only a few entries and for the Lagrange multiplier it holds $ \Vec{\lambda}\in\R^m$ with $m\ll N$.

\subsection{Contact Treatment}
After the general reduction method has been outlined, we introduce a more accurate way of choosing the basis vectors for the reduced space of the displacements. The underlying idea for increasing the accuracy of a reduction approach for contact problems is to partition the nodal variables into, so called, master and slave nodes. One of the most prominent partitioning methods is the Guyan reduction, also called \textit{static condensation}, see \cite{guyan1965reduction}. The slave nodes often possess much smaller local stiffness compared to the master nodes and do not inherit any noticeable change of motion, a fact that is hidden behind this technique \cite{simeon}. The combination of static condensation with a reduction method performed on the slave nodes results into the Craig-Bampton method \cite{craigbampton}, a well-known technique usually applied to contact problems with linear constraints, see \cite{manvelyan2021efficient}. We, however, want to extend this method in a way that will allow us to use the same principle for reducing \eqref{eq:KKT_transient} including the nonlinear term $\Mat{S}(\Vec{\lambda}).$ 
More precisely, let us assume that the displacement variables and the system matrices in \eqref{eq:KKT_transient} are already partitioned, i.e.,
\begin{align}
\Mat{M} = 
\begin{pmatrix}
\Mat{M}_{MM} & \Mat{M}_{MS}\\
\Mat{M}_{SM} & \Mat{M}_{SS}
\end{pmatrix},
\quad
\Mat{K} = 
\begin{pmatrix}
\Mat{K}_{MM} & \Mat{K}_{MS}\\
\Mat{K}_{SM} & \Mat{K}_{SS}
\end{pmatrix},
\quad
\Mat{f} =\big( \Mat{f}_M, \Mat{f}_S \big)^T.
\end{align} 
Since we want to preserve the constraints during the reduction, a similar partitioning is performed for the constraint matrices
\begin{align}\label{const_mat}
\Mat{D}_k = 
\begin{pmatrix}
\Mat{D}_{k,MM} & \Mat{D}_{k,MS}\\
\Mat{D}_{k,SM} & \Mat{D}_{k,SS}
\end{pmatrix},
\quad
\Mat{c}_k = \big(\Mat{c}_{k,M}, \Mat{c}_{k,S}\big)^T,
\quad k=1,\dots,m.
\end{align}
In Section \ref{sec:contact}, the constraint matrices \eqref{const_mat} were derived. We can observe that the constraints apply only to the master nodes, implying 
\begin{align}\label{matzero1}
&\Mat{D}_{k,MS} = \Mat{0},\
\Mat{D}_{k,SM} = \Mat{0},\\ 
&\Mat{D}_{k,SS} = \Mat{0},\
\Mat{c}_{k,S} = \Mat{0}\label{matzero2}
\end{align}
for all $k = 1, \ldots, m.$
Therefore, it follows that
\begin{align}\label{slave_stiff}
\Mat{S}_{SS}(\Vec{\lambda}) = \Mat{K}_{SS} - \sum_{k=1}^{m} \lambda_{k}(\Mat{D}_{k,SS} + \Mat{D}_{k,SS}^T) = \Mat{K}_{SS}.
\end{align}
The equation \eqref{slave_stiff} combined with the vanishing block-matrices corresponding the slave nodes indicate that, as expected, the dynamics of the slave nodal variables are not directly constrained. Even more, if we keep the configuration of the contact nodes unchanged, i.e, $ \Vec{q}_M = 0, $  the dynamic motion of the slave nodes is then given by
\begin{equation} \label{eq:slave}
\Mat{M}_{SS}\ddot{\Vec{q}}_{S} + \Mat{K}_{SS}\Vec{q}_S = \Mat{f}_S .
\end{equation}
Instead of reducing the full system, we want to reduce only the slave nodal variables. Therefore, we compute the transformation matrix of the slave system \eqref{eq:slave} using the Arnoldi method, which yields
\begin{equation}\label{eq:Bkrylov_S}
\Mat{Q}_S:=\left\{\Mat{K}_{SS}^{-1}\Vec{f}_S, 
\Mat{K}_{SS}^{-1}\Mat{M}_{SS} \Mat{K}_{SS}^{-1}\Vec{F}_S,\cdots,
(\Mat{K}_{SS}^{-1}\Mat{M}_{SS})^{n-1} \Mat{K}_{SS}^{-1}\Vec{f}_S \right\}.
\end{equation}
Now we can construct a smaller space approximating the displacement variables outside the contact area. Nevertheless we still need a coupling condition between the slave and the master variables. In order to achieve this we assume that the impact of the acting force on the slave nodes is negligible, implying a zero acceleration for the slave nodes. This assumption leads to a static condensation between the slave and the master nodes, which due to \eqref{slave_stiff} reads
\begin{equation} \label{eq:coupling}
\Mat{K}_{SM} \Vec{q}_{M} + \Mat{K}_{SS} \Vec{q}_{S} = 0.
\end{equation}

Using \eqref{eq:coupling}, we obtain the complete transformation matrix 
that includes the coupling term of the master and slave nodes,
\begin{align}\label{eq:defQCB}
\Mat{Q}_{CB} =
\begin{pmatrix}
\Mat{I}_M &  0 \\
-\Mat{K}_{SS}^{-1}\Mat{K}_{SM} & \Mat{Q}_S
\end{pmatrix}.
\end{align}
The partitioning of the master and slave nodes in the fashion of the Craig-Bampton method provides an alternative way
to compute the reduced displacements in \eqref{eq:system_time_red}. Although the contact condition is quadratic, the static condensation remains unaffected and does not depend on the Lagrange multiplier, adopting the linear form \eqref{eq:coupling}. On the other hand, since the contact nodes are preserved, the individual impact of the contact elements is comprised within the reduced model. A similar result for contact problems with linear constraint is obtained in \cite{manvelyan2021efficient}. 
Regardless whether we employ $\Mat{Q}$ from the reduction (\ref{eq:Bkrylov}) or 
$ \Mat{Q}_{CB} $, the computation of the transformation can be done in an a priori offline phase and does not require solution trajectories, i.e., snapshots, of the full system.

\subsection{The reduction scheme}\label{reduction_scheme}
After these preparations, we present the final reduction scheme. Note that all the reduced matrices can be computed beforehand in an offline procedure. Once the transformation matrix \eqref{eq:defQCB} is computed the full matrices can be reduced via
\begin{align}
&\hat{\Mat{M}} =\Mat{Q}_{CB}^T \Mat{M}\Mat{Q}_{CB} \in \R^{n\times n}, \quad \hat{\Mat{K}}=\Mat{Q}_{CB}^T \Mat{K}\Mat{Q}_{CB} \in \R^{n\times n}, \label{red_MK}\\
&\hat{\Mat{S}}(\Vec{\lambda}) = \Mat{Q}_{CB}^T\Mat{S}(\Vec{\lambda})\Mat{Q}_{CB} =\hat{\Mat{K}}  - \sum_{k=1}^{m} \lambda_{k} (\hat{\Mat{D}}_k + \hat{\Mat{D}}_k^T)  \in \R^{n\times n}, \label{red_S}
\end{align}
whereas the matrices $\hat{\Mat{D}}_k$ comprise solely the constraint  block-matrices corresponding to the master nodes appended by $\text{rank}(\Mat{Q}_{SS})$-many zero columns and rows, i.e.,
\begin{align}
\hat{\Mat{D}}_k = \Mat{Q}_{CB}^T \Mat{D}_k\Mat{Q}_{CB} =
\begin{pmatrix}
\Mat{D}_{k,MM} & \Mat{0}\\
\Mat{0} &\Mat{0}
\end{pmatrix} \in \R^{n \times n} \quad \text{for } k = 1, \dots, m.
\end{align}
The same argument applies to the linear term of the constraint \eqref{const1},
\begin{align}
\hat{\Vec{c}}_k = \Mat{Q}_{CB}^T\Vec{c}_k = \big(\Vec{c}_{k,MM}, \Vec{0}\big)^T \in \R^n,\quad k = 1, \ldots, m. \label{red_c}
\end{align}
The Arnoldi method requires only the position vector of the right hand-side and the rearrangement of the position vector according to the same partitioning results in
\begin{align}\label{red_f}
\hat{\Vec{f}}_{\text{pos}} = \Mat{Q}_{CB}^T\Vec{f}_{\text{pos}}\Mat{Q}_{CB} = \big(\Vec{f}_{\text{pos}, M}, \hat{\Vec{f}}_{\text{pos}, S}\big)^T \in \R^n.
\end{align}
Note that once the matrices $\hat{\Mat{K}}$ and  $\hat{\Mat{D}}_k$ are computed, the stiffness term $\hat{\Mat{S}}(\Vec{\lambda})$ can be recovered be means of \eqref{red_S} for each time step during the online computation.\\ 
Based on the notations with respect to the reduced system introduced so far, the reduced  displacement vector can be formulated as
\begin{equation}\label{displacement_red}
\hat{\Vec{w}}_{t_{i+1}}(\Vec{\lambda}_{t_{i+1}}) = (\hat{\Mat{M}}+h^2\hat{\Mat{S}}(\Vec{\lambda}_{t_{i+1}}))^{-1}(h^2\hat{\Vec{f}}_{t_{i+1}} + h^2\hat{\Mat{C}}^T\Vec{\lambda}_{t_{i+1}} +2\hat{\Mat{M}}\hat{\Vec{w}}_{t_{i}} - \hat{\Mat{M}}\hat{\Vec{w}}_{t_{i-1}}).
\end{equation}
Finally, the constraint function of the reduced space $\hat{\Vec{F}}:\mathbb{R}^m \to \mathbb{R}^m$ reads
\begin{align}\label{F_temp_red} 
	\hat{\Vec{F}}_k(\Vec{\lambda}_{t_{i+1}}) &=\hat{\Vec{w}}_{t_{i+1}}(\Vec{\lambda}_{t_{i+1}})^T\hat{\Mat{D}}_k \hat{\Vec{w}}_{t_{i+1}}(\Vec{\lambda}_{t_{i+1}}) + \hat{\Vec{c}}_k^T \hat{\Vec{w}}_{t_{i+1}}(\Vec{\lambda}_{t_{i+1}}) + b_k, \quad k = 1,\ldots,m.
\end{align}

Similarly to the case of the full order model, the NCP formulation for the constraint function $\hat{\Vec{F}}$ of the reduced system can be derived by inserting the reduced system matrices \eqref{red_MK}-\eqref{red_c} into \eqref{F_temp}. The approximating LCP problem can be stated by means of 
\begin{equation}\label{eq:LCPmatrices}
\begin{aligned}
\hat{\Mat{A}}^{l-1} &=\nabla \hat{\Vec{F}}(\Vec{z}_{l-1}),\\
\hat{\Mat{B}}^{l-1} &= \hat{\Vec{F}} (\Vec{z}_{l-1}) - \nabla \hat{\Vec{F}}(\Vec{z}_{l-1})\Vec{z}_{l-1},
\end{aligned}
\end{equation}
which are computed during the linearization procedure of \eqref{F_temp_red}. Similar to the full model case, we set $\Vec{z}_l = \Vec{\lambda}_{t_{i+1}}$, implying that the Lagrange multiplier, unlike the displacement variable, remains in the full dimensional form. The latter holds due to the fact that with our reduction approach we preserve the contact shape. The final form of the LCP problem at the $l$-th iteration for a fixed time step reads
\begin{equation} \label{eq:solveLCP3}
	\begin{array}{rcc}
	\hat{\Vec{B}}^{l-1} + \hat{\Mat{A}}^{l-1} \Vec{z}_{l} & \geq & 0, \\
	\Vec{z}_{l} & \geq & 0, \\
	\Vec{z}_{l}^T (\hat{\Vec{B}}^{l-1} + \hat{\Mat{A}}^{l-1} \Vec{z}_{l}) & = & 0.
	\end{array}
\end{equation}
 The system \eqref{eq:solveLCP3} represents the linearized adjoint problem of the reduced system \eqref{eq:reduced_equation}. In the terminology of differential-algebraic equations,  
$\Vec{q}$ stands for the differential variables while  
$\Vec{\lambda}$ is the algebraic variable. Based on the numerical analysis performed so far, we can deduce that there is no practical dependency on the starting point $\lambda_0.$ However, further investigations of this theoretic aspect might provide extra insight. In case of a constant Jacobian matrix $\nabla\hat{\Vec{F}}(\Vec{z}_0)$ throughout the fixed-point iterations only linear convergence can be achieved. In order to increase the convergence rate, an update of the Jacobian matrix is necessary. In that case quadratic convergence is guaranteed.
\begin{algorithm}[H]
    \caption{\label{alg:iterLCP} Reduction of the quadratic contact problem \eqref{eq:KKT_transient}}
    \begin{algorithmic}[1]
     \State \mbox{Input:} The system matrices and vectors $\Mat{M}, \Mat{K}, \Vec{f}_{\text{pos}}, \Mat{D}_k, \Vec{c}_k, b_k,\ k = 1,\ldots,m,\ TOL > 0 $
	 \State \mbox{Output:} $\Vec{q}_{t_i}$ and $\Vec{\lambda}_{t_i}$  for $i \in \{1,\ldots,n_T\}$
	 \State Perform partitioning between master and slave nodes.
	 \State Compute the reduction matrix $\Mat{Q}_{CB}$ by \eqref{eq:defQCB}.
	 \State Compute the reduced system by \eqref{red_MK}-\eqref{red_f}.
	 \State Initialization of $\hat{\Vec{w}}_{t_{i}}(\Vec{\lambda}_{t_{i}})$ and $\Vec{\lambda}_{t_i}$ for $i = 1,2.$
	 \For{$i=2,\dots, n_T$}
	    \State $\Vec{z}_0 = \Vec{\lambda}_{t_{i-1}}$
		\For{$l=1,2,\dots$}
		\State Compute $\hat{\Mat{A}}^{l-1}$ and $\hat{\Mat{B}}^{l-1}$ by \eqref{eq:LCPmatrices},
		\State Compute $\Vec{z}_{l}$ as a solution of \eqref{eq:solveLCP3} by Lemke's method.
		\If{$\|\Vec{z}_{l}-\Vec{z}_{l-1}\| < TOL$}
		\State $l^* = l$
		\State \textbf{break}
		\EndIf
		\EndFor
		\State $\Vec{\lambda}_{t_i} = \Vec{z}^{l^*}$
		\State Compute the reduced displacement vector $\hat{\Vec{w}}_{t_{i}}(\Vec{\lambda}_{t_{i}})$ by \eqref{displacement_red}.
	\EndFor
	\State Compute $\Vec{q}_{t_i} = \Mat{Q}_{CB}\hat{\Vec{w}}_{t_{i}}(\Vec{\lambda}_{t_{i}}),\ i = 1,\ldots n_T.$
	 \end{algorithmic}	
\end{algorithm}

We briefly summarize here the properties of the Algorithm \ref{alg:iterLCP}. The reduced model is computed once in an offline phase and it can be utilized for arbitrary right-hand sides. The contact shape is precomputed in the offline phase as well, since it depends on the acting force. The performance of the reduction method depends on the size of the contact zone $m$. The latter holds due to the fact, that the dual problem of the size $m$ is solved within several iterations in each time step. However, according to the theory \cite{josephy1979newton,pang1982iterative}, quadratic convergence is guaranteed for the fixed-point iterations in Algorithm \ref{alg:iterLCP}. Therefore, the convergence criterion can be achieved after a few iterations. Based on these observations and facts, the Algorithm \ref{alg:iterLCP}
can be expected to be highly efficient for contact problems with a small contact zone compared to the overall structure.

\section{Applications}
\label{sec:applications}
In this section we focus on the application of the new reduction scheme and provide two computational examples: The first example is a plane stress self-contact problem serving solely as a demonstrative example without any practical relevance. The second example represents a train wheel-rail contact problem with a potential real-life application in the railway industry. In both scenarios the sliding friction at the contact interface is neglected. The LCPs within the reduction scheme are solved by the Lemke's method taken from the open source library \cite{lemke}.

\subsection{Simple Crack in a Square}
First, we apply our approach to a two-dimensional plane unit square  $\Omega = [0,1] \times [0,1]$ with a self-contact under the assumption of plane stress. For the spatial  discretization  standard bilinear shape functions are used on quadrilateral elements.
The finite element method  and the reduction scheme are implemented within a custom-written Python script. The material parameters $\rho =1,\  E = 1000,\ \nu = 0.3 $ are given in dimensionless form. We consider a zero Dirichlet boundary condition on the left edge of the plane square. The square is discretized by 1600 quadrilateral elements, which in total lead to $N = 3386$ DOF. The interior tear defines the contact zone and is given by a fixed number $ m = 25 $ of discrete gird points. A special care is taken for the data-structure of the nodes at the contact interface. We make use of the node-to-segment technique \cite{wriggers2004computational}. In particular, the left edge is composed of a fixed number of segments. Each of such segments is predefined by its start and end nodes. In contrast, the right side of the tear is given by grid points that coincide with the finite element nodes. 

\begin{figure}[H]
\centering
\begin{subfigure}[b]{.5\textwidth}
\includegraphics[width=\textwidth]{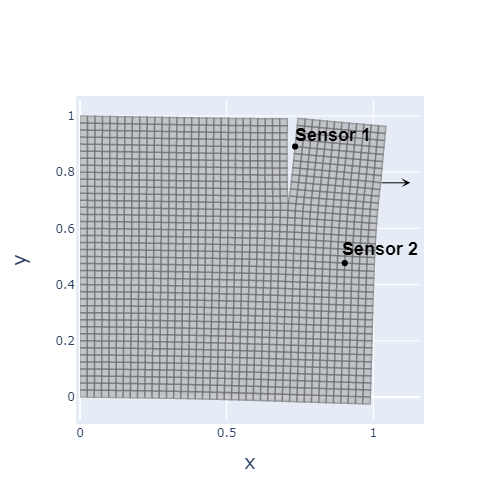}
\caption{No contact at the crack\label{square1}}
\end{subfigure}%
\begin{subfigure}[b]{.5\textwidth}
\includegraphics[width=\textwidth]{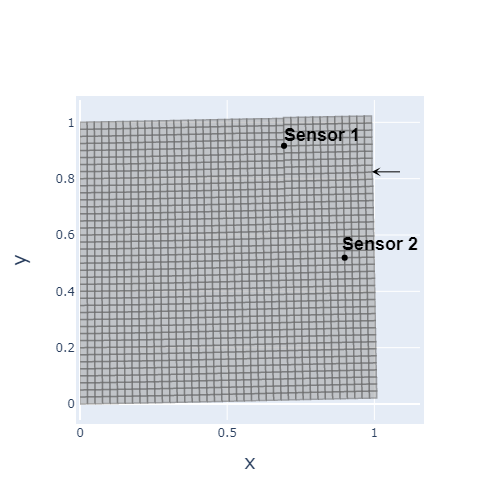}
\caption{Sliding contact at the crack\label{square2}}
\end{subfigure}
	\caption{Two scenarios for the plane square.}\label{fig:square}
\end{figure}

This kind of data-structure allows us to formulate a non-penetration condition during a local sliding movement at the tear of the plane square. For this example, there is a unique assignment between the contact nodes and the contact segments. Note, that at the the initial time point the tear is closed i.e., there is no gap between the contact nodes and the contact segments.\\
The right side of the domain is stressed by a horizontal nodal oscillating force: 
\begin{equation}\label{load1}
\Vec{f}(t) = \left( 1.5\sin(0.1\pi t), 0 \right)^T \in \R^2, \quad t \in [0,1].
\end{equation}
The application of \eqref{load1} has an impact on the contact interface leading to an opening and closing the tear over time. Whenever the tear is closed, a sliding motion occurs at the contact zone.
The problem setup and the mesh are illustrated in Fig. \ref{fig:square}. In order to compare the displacements of the full and the reduced model, two sensors are attached to the square. One of those sensors is placed on the contact interface allowing us to track the contact pressure and the gap function.\\
We eliminate the fixed degrees of freedom a priori before starting with the reduction method. In total, the reduced model has $n = 53$ DOF, which results from the combination of $ 50 $ contact  displacement nodes and $3$ additional slave variables.

\begin{figure}[H]
\centering
\begin{subfigure}[b]{.5\textwidth}
\includegraphics[width=\textwidth]{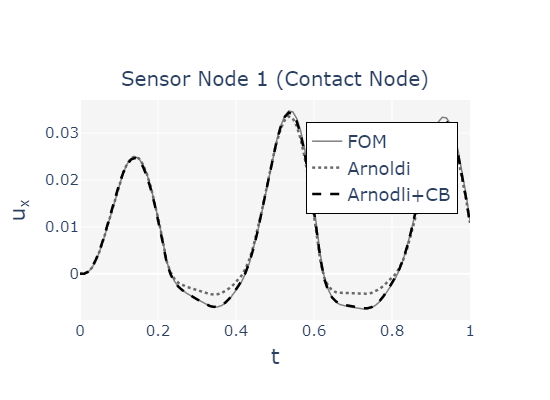}
\caption{\label{node1_x}{The x-displacements}}
\end{subfigure}%
\begin{subfigure}[b]{.5\textwidth}
\includegraphics[width=\textwidth]{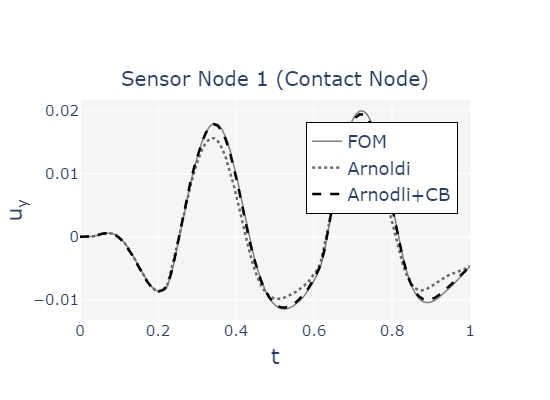}
\caption{\label{node1_y}{The y-displacements}}
\end{subfigure}
	\caption{Comparison of FOM and two ROMs (with and without contact treatment).}\label{fig:square_node1}
\end{figure}

\begin{figure}[H]
\centering
\begin{subfigure}[b]{.5\textwidth}
\includegraphics[width=\textwidth]{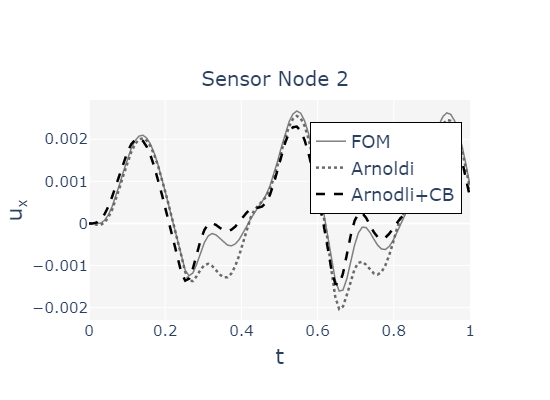}
\caption{\label{node2_x}{The x-displacements}}
\end{subfigure}%
\begin{subfigure}[b]{.5\textwidth}
\includegraphics[width=\textwidth]{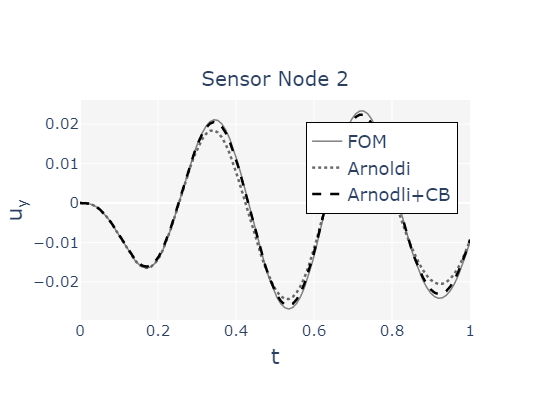}
\caption{\label{node2_y}{The y-displacements}}
\end{subfigure}
	\caption{Comparison of FOM and two ROMs (with and without contact treatment). }\label{fig:square_node2}
\end{figure}

First, the reduction approach is performed by means of an Arnoli reduction matrix without any contact treatment. Afterwards, we compute the displacements resulting from the reduction method combined with the Craig-Bampton splitting technique. For the sake of a fair comparison, we choose the same number of basis functions in both cases. When we omit the splitting, the basis consists solely of Krylov basis vectors, whereas in case of splitting, the basis is extended by the contact nodes. Similar to \cite{manvelyan2021efficient}, the advantages of additional contact treatment can be clearly observed. Both solutions resulting from the corresponding reduction approach are compared with the full displacement variables in Fig.~\ref{fig:square_node1} and in Fig.~\ref{fig:square_node2}.

\begin{figure}[H]
\centering
\begin{subfigure}[b]{.5\textwidth}
\includegraphics[width=\textwidth]{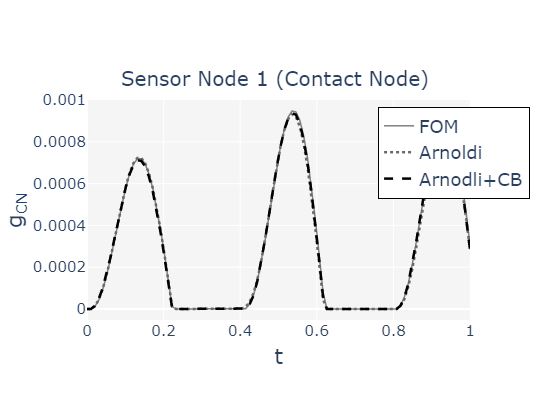}
\caption{\label{cb_node2_x}{The gap function}}
\end{subfigure}%
\begin{subfigure}[b]{.5\textwidth}
\includegraphics[width=\textwidth]{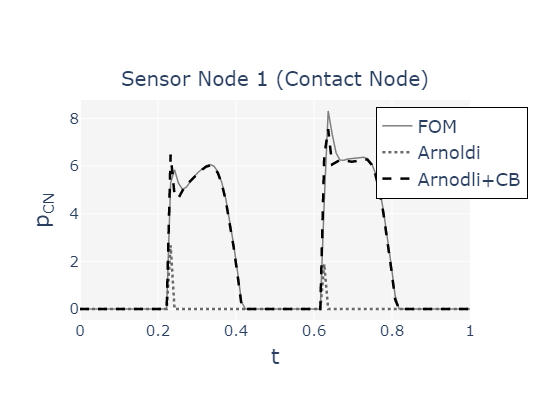}
\caption{\label{cb_node2_x_lmb}{The contact pressure}}
\end{subfigure}
	\caption{Complementarity between the gap funtion and the contact pressure.}\label{fig:gap_vs_pressure}
\end{figure}

Furthermore, the contact pressure and its counterpart, the gap function, deserve attention. We denote $g_{CN} = \Vec{u}^T \Mat{D}_{\hat{k}} \Vec{u} + \Vec{c}_{\hat{k}}^T \Vec{u} + b_{\hat{k}}$ and
$p_{CN} = \lambda_{\hat{k}}$, where ${\hat{k}}$ denotes the sensor node at the contact zone. We observe in Fig.~\ref{fig:gap_vs_pressure}, that, in particular, when computed by the ROM without splitting, the pressure is not recovered, whereas the ROM with a splitting provides the contact pressure, that agrees well with the FOM. Both approaches guarantee positivity of the contact pressure whenever the contact is activated. Furthermore, the switching between an active contact with a nonzero pressure and positive gap with vanishing pressure can be nicely observed, i.e., it holds $g_{CN}p_{CN} = 0$ for all time steps.
\begin{figure}[H]
\centering
\includegraphics[width=\textwidth]{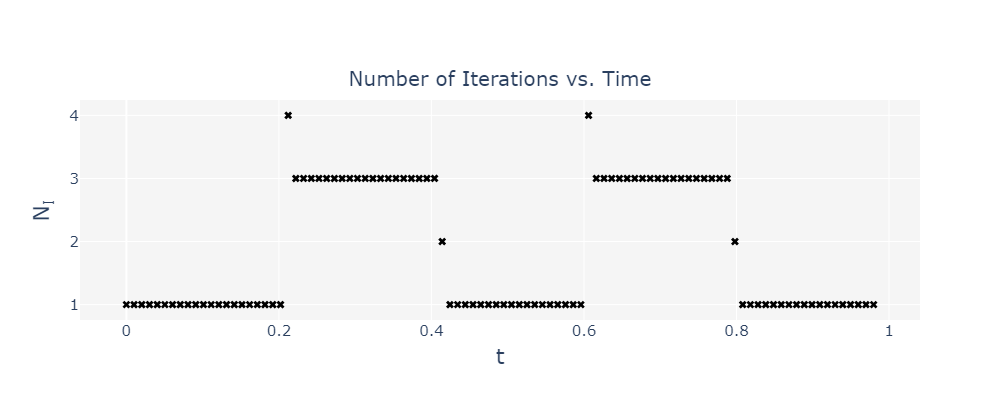}
\caption{The number of the fixed-point iterations for solving a NCP problem over time.}\label{fig:iterations}
\end{figure}
Overall, the main reason for this behaviour is that by using only Arnoldi for computation of the transformation matrix, the reduced space does not comprise the dominant impact of the contact nodes. Moreover, our reduction approach is highly efficient in case of a small contact zone compared to the total volume and since the reduced space is at least as large as the contact area, we do not loose the computational savings by adding the contact splitting to our reduction approach.\\
Finally the convergence behaviour of the NCPs is addressed. When the contact is inactive, only one fixed-point iteration of LCP suffices to solve the non-penetration condition. However, in case of an active contact condition, 3 or 4 iterations are needed for each time-step, see Fig.~\ref{fig:iterations}. This observations agree well with the theory implying quadratic convergence for solving NCPs during each time step.

\subsection{Wheel-Rail Contact Problem}
\begin{figure}[H]
\centering
\begin{subfigure}[b]{.2\textwidth}
\includegraphics[width=\textwidth]{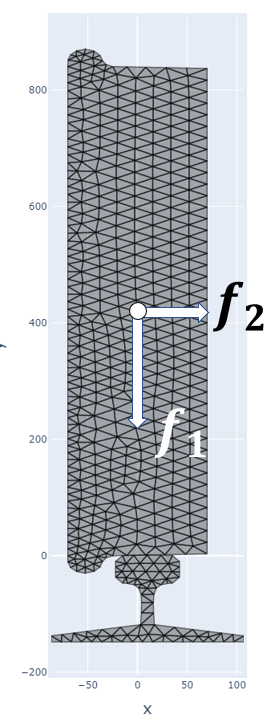}
\caption{Full model \label{wheelrail1}}
\end{subfigure}%
\begin{subfigure}[b]{.3\textwidth}
\includegraphics[width=\textwidth]{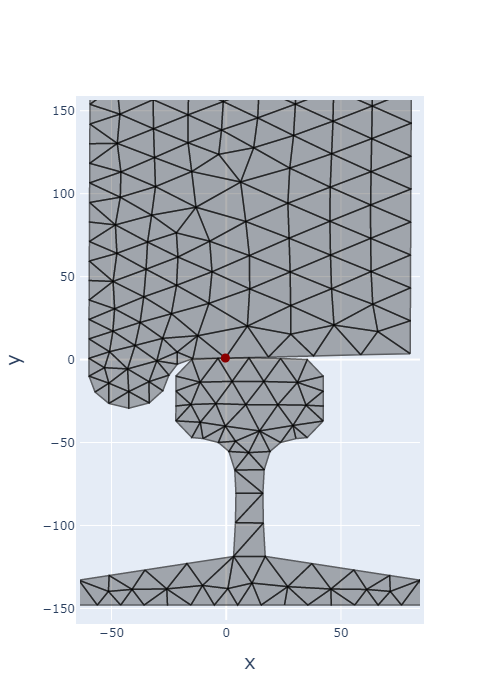}
\caption{Contact area \label{sensor}}
\end{subfigure}
	\caption{A two-dimensional cross section of the wheel-rail contact}
\end{figure}

To conclude the computational section we present the train wheel-rail contact problem, a prominent example of real-life application in contact mechanics. Both the wheel and the rail are composed of alloy steel and discretized within the FEM-framework of NX Simcenter 12.0. Within the scope of this work we solely consider the two-dimensional cross section of the wheel-rail.  Furthermore, the rail is fixed on the bottom over time and the rotational degrees of freedom are neglected. Both the wheel and the rail share the same Young modulus, $E = 206940$ MPa, and the Poisson coefficient $ \nu = 0.288.$ Two outer loads are applied to the wheel: The first one is the superposition of all vertical forces pointing towards ground (gravity, train mass, etc.). The second force represents the centripetal force defined by a frequency and load magnitude (Fig. \ref{wheelrail1}). The force $\Vec{f}_2$ arises naturally during the train ride and the corresponding frequency usually depends on the train velocity. All trajectories are computed at the sensor node placed on the contact zone of the wheel, see Fig. \ref{sensor}.

In contrast to the previous example, in this case the underlying dynamical system represents a two-body problem. The mass and stiffness matrices consist of two diagonal blocks where the upper diagonal block refers to the wheel whereas the lower diagonal block refers to the rail body. The coupling between the two bodies is given by the contact condition \eqref{const1}-\eqref{const3}.\\

\begin{figure}[H]
\centering
\begin{subfigure}[b]{.5\textwidth}
\includegraphics[width=\textwidth]{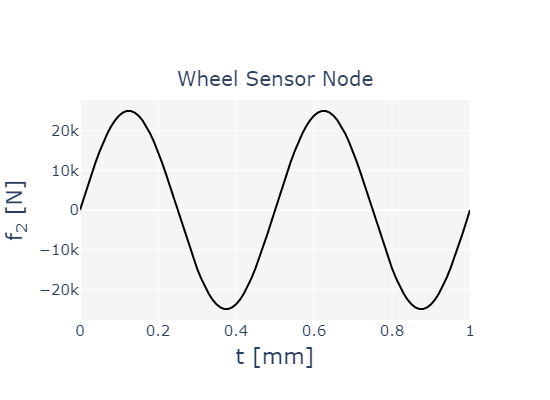}
\caption{\label{force}{$\mathit{f} = 4.0$ Hz, $f_{mag} = 25000 $ N}}
\end{subfigure}%
\begin{subfigure}[b]{.5\textwidth}
\includegraphics[width=\textwidth]{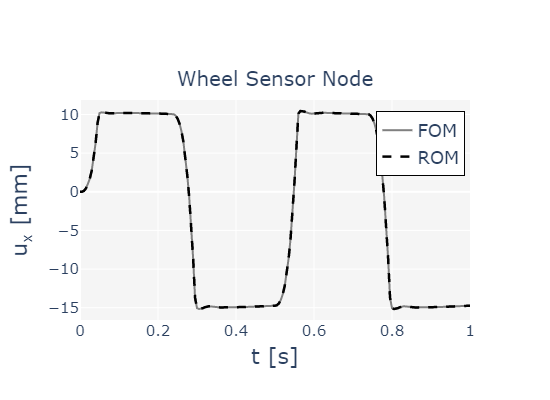}
\caption{\label{ux}{The x-displacements}}
\end{subfigure}
	\caption{The outer horizontal force and the displacements}\label{fig:rw1}
\end{figure}
\begin{figure}[H]
\centering
\begin{subfigure}[b]{.5\textwidth}
\includegraphics[width=\textwidth]{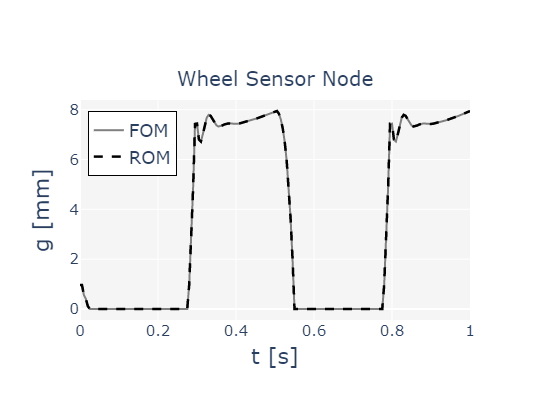}
\caption{\label{gap}{}The gap function}
\end{subfigure}%
\begin{subfigure}[b]{.5\textwidth}
\includegraphics[width=\textwidth]{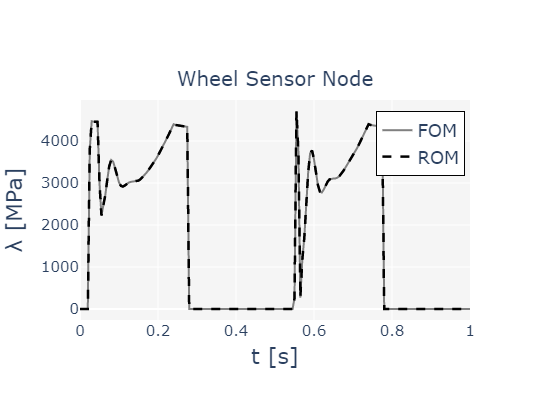}
\caption{\label{lambda}{}The contact pressure}
\end{subfigure}
	\caption{Complementarity between the gap function and the contact pressure.}\label{fig:rw2}
\end{figure}

In Fig.~\ref{fig:rw1}, the x-displacements of the full and the reduced model, that are computed at the sensor node, are depicted. The comparison of the both trajectories shows a very good agreement between the reduced and the full model. Moreover, the contact pressure and the gap function calculated at the same nodes are shown in Fig.~\ref{fig:rw2}. Note, that between the wheel and the rail a periodic sliding movement occurs. Moreover, the complementarity between the gap function and the contact pressure retains over time.

\begin{figure}[H]
\centering
\begin{subfigure}[b]{.5\textwidth}
\includegraphics[width=\textwidth]{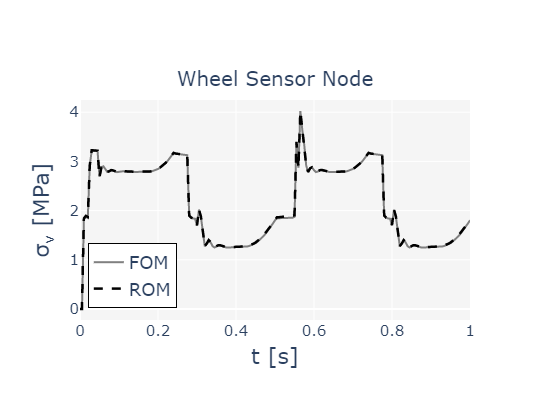}
\caption{\label{mises}{von Mises stress}}
\end{subfigure}%
\begin{subfigure}[b]{.5\textwidth}
\includegraphics[width=\textwidth]{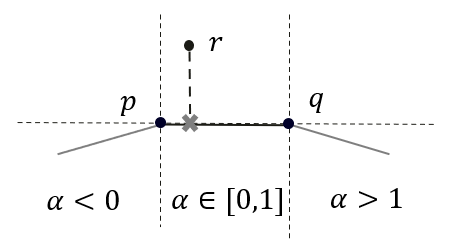}
\caption{\label{update}{Contact update technique}}
\end{subfigure}
	\caption{The nodal stress and the illustration of the contact update technique.}\label{fig:rw3}
\end{figure}

Due to the additional contact treatment by Craig-Bampton within our reduction framework, a specific form of contact update is possible. Since the contact nodes and the contact segments are recovered by the reduced model at each time step, the node-segment correspondence can be reestablished, if necessary. Various suitable criterions can be used for such an update. In particular, in case of sliding, nodes can leave the surrounding area of the corresponding segment over time. Therefore, the main idea is based on finding the nearest segment of each contact node.\\
In Algorithm \ref{alg:contactUpdate} one iteration of a contact update that can be integrated into our reduction algorithm is introduced. Within this algorithm we make use of the following geometrical observation:
For the projection of the node $r \in \R^2$ to the line going trough the corresponding segment nodes $p, q \in \mathbb{R}^2$ it holds
\begin{align}
\Pi(r) = (1-\alpha)p +  \alpha q,\quad \alpha \in \R.
\end{align}
Depending on the value of $\alpha$ we have tree different scenarios for the node $r$, see Fig.~\ref{update}. As soon as $\alpha$ leaves the interval $[0,1],$ an update step either to the right ($\alpha >1$) or to the left segment ($\alpha < 0$) is preformed. However, a certain tolerance is added for smoothing the discontinuity arising due the update technique, see Algorithm \ref{alg:contactUpdate}. Note, that segments are stored within a linear list during the offline procedure. Moreover, there is a linear selecting function, that maps the contact nodes to the contact segments. The latter can be updated after each time step and is the output of Algorithm \ref{alg:contactUpdate}.

\begin{algorithm}[H]
    \caption{\label{alg:contactUpdate} Contact Update for \eqref{eq:KKT_transient}}
    \begin{algorithmic}[1]
     \State \mbox{Input:} List of segments $\mathcal{S} = (s_1, s_2,\ldots, s_{k_{\mathcal{S}}}),$ list of nodes $\mathcal{N} = (p_1,p_2\ldots p_{k_{\mathcal{N}}})$, 
     \State \mbox{Input:} The old selecting function $\Phi:\mathcal{N} \to \mathcal{S},\ TOL \in (0,1)$
	 \State \mbox{Output:} The new selecting function $\tilde{\Phi} : \mathcal{N} \to \mathcal{S}$
	 \For{$i=1,\dots k_{\mathcal{N}}$}:
	    \State Get $s_j = \Phi(p_i).$
		\State Compute the projection  $\Pi(p_i)$ w.r.t. the segment $s_j = (p^1_{s_j},p^2_{s_j}).$
		\State Solve $\alpha$ for $ \Pi(p_i)  = (1-\alpha)p^1_{s_j} + \alpha p^2_{s_j}. $
		\If{$\alpha <- TOL$ and $j>1$}
		    \State Set $\tilde{\Phi}(p_i) = s_{j-1}$
		\EndIf
		\If{$\alpha > 1+ TOL$ and $j<k_{\mathcal{S}}$}
		    \State Set $\tilde{\Phi}(p_i) = s_{j+1}$
		\EndIf
	\EndFor
	 \end{algorithmic}	
\end{algorithm}

Finally, we want to compute the stress at the sensor node both for the full as well as for the reduced model. From the definition of stress \eqref{stress_strain}, it follows that the latter is a function of the gradient of the displacements. Within a FEM-framework the gradient of the displacement function $\nabla\Vec{u}_{i}$ at one node $\Vec{p}_i$ can be approximated by the averaged sum of the gradients computed on all elements $e_j$, that are adjacent to $\Vec{p}_i,$ i.e.,
\begin{align}
\nabla \Vec{u}_{i} = \frac{1}{m}\sum_{j = 1,\ldots, m}\nabla  \Vec{u}_{|{e_j}}(\Vec{p}_i),
\end{align}
where  $\Vec{u}_{|{e_j}}$ denotes the displacement function restricted on the element $e_j.$ In case of triangle linear elements straightforward calculations show that for an element $e_j$ with the nodes $\Vec{p}^j_{i_1},\Vec{p}^j_{i_2}, \Vec{p}^j_{i_3}$ it holds
\begin{align}
\nabla \Vec{u}_{|{e_j}} = \Mat{D}_j\Mat{B}_j(\Vec{p}^j_{i_1},\Vec{p}^j_{i_2}, \Vec{p}^j_{i_3}),
\end{align}
where $\Mat{D}_j =\big(\Vec{u}^j_{i_1},\Vec{u}^j_{i_2}, \Vec{u}^j_{i_3}\big) \in \R^{2\times 3}$ contains the nodal displacements of the nodes $\Vec{p}^j_{i_1},\Vec{p}^j_{i_2}, \Vec{p}^j_{i_3}$ and $\Mat{B} = \Mat{B}(\Vec{p}^j_{i_1},\Vec{p}^j_{i_2}, \Vec{p}^j_{i_3})  \in \R^{3\times 2}$ depends solely on the initial coordinates of the same nodes. The entries of the matrix $\Mat{B}$ can be easily calculated and are omitted here.
In summary, in order to compute the stress at a node $\Vec{p_i}$, the coordinates and the displacements of all direct neighbouring nodes $\Vec{p}_{k_1},\Vec{p}_{k_2}, \Vec{p}_{k_{m_i}}$ are required. In particular, in our example for the sensor node at the contact zone, there are $m_i = 3$ neighbour nodes. Finally, there are two possible ways to compute the stress from the reduced displacements. One way is to track the rows $i^x,\ i^y$ of the reduction matrix $\Mat{Q}_{CB}$ defined in \eqref{eq:defQCB} fulfilling
\begin{align}
    \Vec{u}_i = \big(\Mat{Q}_{CB}\big)_{i^x, i^y}\Vec{q},
\end{align}
where $\Vec{q} \in \R^N$ is the full-dimensional discretized displacement vector and the indices $i^x,\ i^y$ stand for the x- and y-coordinates of the node $\Vec{p}_i.$ In this case in total $2(m_i+1)$-many rows of the matrix $Q_{CB}$ has to be tracked. Another way is to make use of the Craig-Bampton technique by adding the three neighbouring nodes to the set of the master nodes. The latter method allows a direct computation of the stress within the reduced space without any need for tracking the transformation matrix $Q_{CB}.$

\section{Conclusions}
\label{sec:conclusions}
In \cite{manvelyan2021efficient} we have presented a novel reduction algorithm for structural mechanical problems with linear constraints. The algorithm is purely physics-based and therefore there is no need for snapshots of the state trajectories. The main idea behind this approach was to apply adjoint methods to the reduced system leading to a LCP, that has to be solved in each time step. Furthermore, special care was taken for preserving the contact shape during the reduction process. We have separated the contact nodes from the remaining nodes by applying the Craig-Bampton contact treatment technique. Overall, the reduction method works especially well for problems, where the number of inner volume nodes is larger compared to the contact nodes.

The main drawback of the method described in \cite{manvelyan2021efficient} is that only node-to-node problems are considered, allowing to address only linear constraints. In the current paper we have generalized the reduction algorithm for contact problems, where also sliding movements can be considered. This reduction scheme is designed for node-to-segment contact conditions and all the benefits of the reduction method for node-to-node problems carry over. Moreover, the node-to segment contact formulation leads to quadratic contact conditions, which in terms of variational formulation are appended to the energy functional by Lagrangian multipliers. Because the Kuhn-Tacker condition remains linear in the state, an explicit expression for the adjoint equations may be derived. Here, the Lagrangian multipliers must fulfill the nonlinear dual system, which represent a NCP and has to be solved for each time step. Furthermore, our computational results show, that the corresponding Newton-like algorithm converges in about three iteration. Therefore, the computational effort of the reduction method during a single time step is still comparatively low. 

The performance of the reduction scheme is demonstrated on two computational examples: The first one demonstrates the sliding capability of the reduced model and allows to
study the NCP-problem in detail.
The second application describes a wheel-rail contact problem. Here, we additionally need a variable mapping from nodes to the segments.

Further research will address friction in the contact condition as well as adapting the corresponding reduction algorithm.

\vspace{1.0cm}
{\bf Acknowledgement}: The authors would like to thank Christoph Heinrich and his analysis group for valuable discussions.

\vspace{1.0cm}
\bibliographystyle{plain}
\bibliography{references}%

\end{document}